\documentclass[10pt,twoside]{article}
\usepackage{amsmath,amsthm}
\usepackage[latin1]{inputenc}
\usepackage{picinpar}
\usepackage{longtable}
\usepackage{amsrefs}
\usepackage[dvips]{graphicx}
\usepackage[mathscr]{eucal}
\usepackage{calc}
\usepackage{ifthen}
\usepackage{dcpic,pictexwd}
\usepackage[all]{xy}
\usepackage{enumerate}
\usepackage{amssymb,latexsym}
\usepackage{amsthm}
\usepackage[dvips]{graphics}
\usepackage{color}
\usepackage{tabularx}
\usepackage{mathtools}
\usepackage{makeidx}
\usepackage{graphicx}
\theoremstyle{definition}
\newtheorem{teor}{Theorem}
\newtheorem{lema}{Lemma}
\newtheorem{prop}{Proposition}
\newtheorem{corol}{Corollary}

\newtheorem{Nota}{Remark}

\begin{document}
\thispagestyle{plain}
\par\bigskip

\begin{center}
\textbf{Dyck Paths Categories And Its Relationships With Cluster Algebras}
\end{center}

\par\bigskip

\begin{centering}
  Agust\'{\i}n Moreno Ca\~{n}adas\\{amorenoca@unal.edu.co}\\
Gabriel Bravo R\'{\i}os \\{gbravor@unal.edu.co}\\
Department of Mathematics\\
Universidad Nacional de Colombia\par\bigskip
\end{centering}

\par\bigskip

\begin{centering}
\textbf{Abstract}\par\bigskip
\end{centering}

\small{Dyck paths categories are introduced as a combinatorial model of the category of representations of quivers of Dynkin type $\mathbb{A}_{n}$. In particular, it is proved that there is a bijection between some Dyck paths and perfect matchings of some snake graphs. The approach allows us to give formulas for cluster variables in cluster algebras of Dynkin type $\mathbb{A}_{n}$ in terms of Dyck paths.}

\par\bigskip

\par\bigskip


\section{Introduction}

In the last few years, researches regarding connections between cluster algebras and different fields of mathematics have been growing. For instance, relationships between cluster algebras, quiver representations, combinatorics and number theory have been reported by Fomin et al., Shiffler et al., K. Baur et al., Assem et al. amongst a great number of mathematicians \cites{Assem, Baur, Canakci, Canakci3, Fomin, Fomin1, Fomin2, Fomin3}.\par\bigskip

Perhaps the Catalan combinatorics (which consists of all the enumeration problems whose solutions are Catalan numbers) is the most appropriated environment for the investigation of cluster algebras of Dynkin type $\mathbb{A}_{n}$. Among all these kinds of problems, it is possible to prove (for example) that the Catalan numbers count \cite{Stanley}:

\begin{enumerate}
\item The number of plane binary trees with $n+1$ endpoints (or $2n+1$) vertices,
\item The number of ways to parenthesize a string of length $n+1$ subject to a non associative binary operation,
\item The number of paths $P$ in the $(x, y)$-plane from $(0,0)$ to $(2n, 0)$ with steps $(1, 1)$ and $(1, -1)$ that never pass below the $x$-axis. Such paths are called \textit{Dyck paths},
\item The number of triangulations of an $(n+3)$ polygon,
\item The number of clusters of a cluster algebra of Dynkin type $\mathbb{A}_{n}$. 
\end{enumerate}

Regarding integer friezes, we point out that Propp in \cite{Prop} reminds that Conway and Coxeter completely classified the frieze patterns whose entries are positive integers, and showed that these frieze patterns constitute a manifestation of the Catalan numbers. Specifically, that there is a natural association between positive integer frieze patterns and triangulations of regular polygons with labelled vertices. According to Baur and Marsh \cite{Baur}, a connection between cluster algebras and frieze patterns was established by Caldero and Chapoton \cite{Caldero1}, which showed that frieze patterns can be obtained from cluster algebras of Dynkin type $\mathbb{A}_{n}$.\par\bigskip

Another example of the use of the Catalan combinatorics as a tool to describe the structure of cluster algebras, was given by Schiffler et al. \cites{Canakci, Canakci3, Musiker}, who found out formulas for cluster variables based on its relations with some triangulated surfaces and perfect matchings of snake graphs. They also proved that there is a way of obtaining the number of perfect matchings of a given snake graph by associating a suitable continued fraction defined by the sign function of the graph.\par\bigskip

Given a non-negative integer $n$ and a triangulation $T$ of a regular polygon with $(n+3)$ vertices. Caldero, Chapoton and Schiffler \cite{Caldero} gave a realization of the category $\mathcal{C}_{C}$ of representations of a quiver $Q_{C}$ associated to a cluster $C$ of a cluster algebra in terms of the diagonals of the $(n+3)$ polygon. They proved that there is a categorical equivalence between the categories $C_{T}$ and $\mathrm{Mod}\hspace{0.1cm}Q_{T}$, where $C_{T}$ is the category whose objects are positive integral linear combinations of positive roots (i.e.,  diagonals that does not belong to the triangulation $T$), whereas $\mathrm{Mod}\hspace{0.1cm}Q_{T}$ denotes the category of modules over the quiver $Q_{T}$ with triangular relations induced by the triangulation $T$.\par\bigskip

Follow the ideas of Caldero, Chapoton and Schiffler, in this paper, a combinatorial model of the category of representations of Dynkin quivers with relations is developed by using Dyck paths. To do that, Dyck paths categories are introduced and it is proved that these categories are equivalent to categories of representations of Dynkin quivers of type $\mathbb{A}_{n}$. This approach allows us to realize perfect matchings of snake graphs as objects of suitable Dyck paths categories, and with this machinery a formula for cluster variables  based on Dyck paths is obtained.\par\bigskip

This paper is distributed as follows; In Section 2, notation and basic definitions to be used throughout the paper are introduced. In Section 3, we define Dyck paths categories and some of its main categorical properties are given in Section 4. In section 5, relationships between objects of the categories of Dyck paths, perfect matchings and cluster algebras are given.

\section{Preliminaries}

\subsection{Cluster Algebras}

 Fomin and Zelevinsky introduced the term cluster algebra in \cite{Fomin}, as a subalgebra of a field of rational functions generated by a set of $n$ cluster variables \cite{Fomin3, Fomin1, Fomin2}. The cluster algebras are in connection with different topics in mathematics, as algebraic combinatorics, Lie theory,  discrete dynamical systems, tropical geometry, and others. Afterwards,  Fomin, Schiffler et al introduced cluster algebras associated to surfaces \cite{Canakci3,Fomin3,Musiker}. 
 \par\bigskip

 For the sake of clarity, we remind here the definition given by Fomin et al. \cite{Fomin1} of a cluster algebra.\par\bigskip

 Let $\mathbb{T}_{n}$ the $n$-regular tree whose edges are labeled by the numbers $1,\dots, n$, so that the $n$-edges incident to each vertex  receive different labels. The symbol $t\stackrel{k}{-}t'$ is used to denote that vertices $t, t'\in\mathbb{T}_{n}$ are joined by an edge labeled by $k$.\par\bigskip

 If $\mathscr{F}$ is a field isomorphic to the field of rational functions over $\mathbb{C}$ (alternatively over $\mathbb{Q}$) in $m$ independent variables, then a \textit{labeled seed of geometric type} over $\mathscr{F}$ is a pair $(\widetilde{x},\widetilde{B})$ where;

 \begin{enumerate}
 \item $\widetilde{x}=(x_{1}, x_{2},\dots, x_{m})$ is an $m$-tuple of elements of $\mathscr{F}$ forming a free generating set; that is; $x_{1}, x_{2},\dots, x_{m}$ are algebraically independent, and $\mathscr{F}=\mathbb{C}(x_{1},\dots, x_{m})$,
 \item $\widetilde{B}=(b_{ij})$ is an $m\times n$ extended skew-symmetrizable integer matrix. $\widetilde{B}$ is said to be the extended exchange matrix of the seed. Its top $n\times n$ submatrix $B$ is the exchange matrix. 
 \end{enumerate}

 Let $(\widetilde{X},\widetilde{B})$ be a labeled seed as above. Take an index $k\in\{1,2,\dots, n\}$. The \textit{seed mutation in direction $k$} transforms $(\widetilde{x}, \widetilde{B})$ into the new labeled seed $\mu_{k}(\widetilde{x}, \widetilde{B})=\widetilde{x}',\widetilde{B}'$ defined as follows;

 \begin{equation}\label{mutation1}
 \begin{split} 
 \widetilde{B}'&=\mu_{k}(\widetilde{B})=(b'_{ij}),
 \end{split}
 \end{equation}

 where
 \[b'_{ij}=
\begin{cases}
-b_{ij},  &   \hspace{0.2cm}\text{if}\hspace{0.1cm}i=k\hspace{0.1cm}\text{or}\hspace{0.1cm} j=k, \\
b_{ij}+b_{ik}b_{k j},  &   \hspace{0.2cm}\text{if}\hspace{0.1cm}b_{ik}>0\hspace{0.1cm}\text{and}\hspace{0.1cm} b_{k j}>0, \\
b_{ij}-b_{ik}b_{k j},  &   \hspace{0.2cm}\text{if}\hspace{0.1cm}b_{ik}<0\hspace{0.1cm}\text{and}\hspace{0.1cm} b_{k j}<0, \\
b_{i j},  &   \hspace{0.2cm}\text{otherwise}. \\
\end{cases}\]

 The \textit{extended cluster} $\widetilde{x}'=(x'_{1},\dots, x'_{m})$ is given by the identifications $x'_{j}=x_{j}$ for $j\neq k$, whereas $x'_{k}\in\mathscr{F}$ is determined by the exchange rule.

 \begin{equation}\label{exchangerule}
 \begin{split}
 x_{k}x'_{k}&=\underset{b_{ik}>0}{\prod}x^{b_{ik}}_{i}+\underset{b_{ik}<0}{\prod}x^{-b_{ik}}_{i}.
 \end{split}
 \end{equation}
 
 A \textit{seed pattern} is defined by assigning a labeled seed $(\widetilde{x}(t),\widetilde{B}(t))$ to every vertex, $t\in\mathbb{T}_{n}$, so that the seeds assigned to the end points of any edge  $t\stackrel{k}{-}t'$ are obtained from each other by the seed mutation in direction $k$. A seed pattern is uniquely determined by one of its seeds.\par\bigskip

 Let $(\widetilde{x}(t),\widetilde{B}(t))_{t\in\mathbb{T}_{n}}$ be a seed pattern as above, and let $\mathscr{X}=\underset{t\in\mathbb{T}_{n}}{\bigcup}x(t)$ be the set of all cluster variables appearing in its seeds. We let the ground ring be $R=\mathbb{C}[x_{n+1},\dots, x_{m}]$ the polynomial ring generated by the \textit{frozen variables}.\par\bigskip

 The \textit{cluster algebra} $\mathscr{A}$ (of geometric type over $R$) associated with the given seed pattern is the $R$-subalgebra of the ambient field $\mathscr{F}$ generated by all cluster variables $\mathscr{A}=R[\mathscr{X}]$.

 \subsubsection{Cluster Algebras From Quivers}

 For quivers, cluster algebras are defined as follows:\par\bigskip

 Fix an integer $n\geq1$. In this case, a seed $(Q, u)$ consists of a finite quiver $Q$ without loops or 2-cycles with vertex set $\{1,\dots, n\}$, whereas $u$ is a free-generating  set $\{u_{1},\dots, u_{n}\}$ of the field $\mathbb{Q}(x_{1},\dots, x_{n})$.\par\bigskip

 Let $(Q, u)$ be a seed and $k$ a vertex of $Q$. The mutation $\mu_{k}(Q, u)$ of $(Q, u)$ at $k$ is the seed $(Q',u')$, where;

 \begin{enumerate}[(a)]
 
\item $Q'$ is obtained from $Q$ as follows;

\begin{enumerate}[(1)]
\item reverse all arrows incident with $k$,
\item for all vertices $i\neq j$ distinct from $k$, modify the number of arrows between $i$ and $j$, in such a way that a system of arrows of the form $(i \stackrel{r}{\longrightarrow} j, i \stackrel{s}{\longrightarrow} k, k \stackrel{t}{\longrightarrow} j) $ is transformed into the system $(i \stackrel{r+st}{\longrightarrow} j, k \stackrel{s}{\longrightarrow} i, j \stackrel{t}{\longrightarrow} k) $. And the system $(i \stackrel{r}{\longrightarrow} j, j\stackrel{t}{\longrightarrow} k, k \stackrel{s}{\longrightarrow} i) $ is transformed into the system $(i \stackrel{r-st}{\longrightarrow} j, i \stackrel{s}{\longrightarrow} k, k \stackrel{t}{\longrightarrow} j) $.  Where, $r$, $s$ and $t$ are non-negative integers, an arrow $i \stackrel{l}{\longrightarrow} j$, with $l\geq0$ means that $l$ arrows go form $i$ to $j$ and an arrow $i \stackrel{l}{\longrightarrow} j$, with $l\leq0$ means that $-l$ arrows go from $j$ to $i$.
\end{enumerate} 

\item $u'$ is obtained form $u$ by replacing the element $u_{k}$ with

\begin{equation}\label{exchangeruleq}
 \begin{split}
 u_{k}&=\frac{1}{u_{k}}\underset{\mathrm{arrows}\hspace{0.1cm}i\rightarrow k}{\prod}u_{i}+\underset{\mathrm{arrows}\hspace{0.1cm}k\rightarrow j}{\prod}u_{j}.
 \end{split}
 \end{equation}

 \end{enumerate}

 If there are no arrows from $i$ with target $k$, the product is taken over the empty set and equals 1. It is not hard to see that $\mu_{k}(\mu_{k}(Q, u))=(Q, u)$.
 In this case the matrix mutation $B'$ has the form

  \[b'_{ij}=
\begin{cases}
-b_{ij},  &   \hspace{0.2cm}\text{if}\hspace{0.1cm}i=k\hspace{0.1cm}\text{or}\hspace{0.1cm} j=k, \\
b_{ij}+sgn(b_{ik})[b_{ik}b_{k j}]_{+},  &   \hspace{0.2cm}\text{else}, 
\end{cases}\]

where $[x]_{+}=\mathrm{max(x,0)}$. Thus, if $Q$ is a finite quiver without loops or 2-cycles with vertex set $\{1,\dots, n\}$, the following interpretations have place:

\begin{enumerate}
\item the clusters with respect to $Q$ are the sets $u$ appearing in seeds, $(Q, u)$ obtained from a initial seed $(Q, x)$ by iterated mutation,
\item the cluster variables for $Q$ are the elements of all clusters,
\item the cluster algebra $\mathscr{A}(Q)$ is the $\mathbb{Q}$-subalgebra of the field $\mathbb{Q}(x_{1},\dots, x_{n})$ generated by all the cluster variables.
\end{enumerate}

 As example, the cluster variables associated to the quiver $Q=1\longrightarrow2$ are:\par\bigskip 

 \begin{centering}
 $\{x_{1}, x_{2}, \frac{1+x_{2}}{x_{1}}, \frac{1+x_{1}+x_{2}}{x_{1}x_{2}}, \frac{1+x_{1}}{x_{2}}\}$.\par\bigskip
 \end{centering}

\subsubsection{Cluster Algebras From Surfaces}

 Let $S$ be a connected oriented $2-$dimensional Riemann surface with nonempty boundary, and let $M$ be a nonempty finite subset of the boundary of $S$, such that each boundary component of $S$ contains at least one point of $M$. The elements of $M$ are called \textit{marked points}. The pair $(S,M)$ is called a  \textit{bordered surface with marked points}. Marked points in the interior of $S$ are called punctures (For technical reasons, we require that $(S,M)$ is not a disk with 1,2 or 3 marked points) \cite{Canakci3}.\\

An \textit{arc} $\gamma$ in $(S,M)$ is a curve in $S$, considered up to isotopy, such that:

\begin{enumerate}
\item [(i)]the endpoints of $\gamma$ are in $M$,
\item [(ii)] $\gamma$ does not cross itself, except that its endpoints, may coincide,
\item [(iii)] except for the endpoints, $\gamma$ is disjoint from the boundary of $S$,
\item [(iv)] $\gamma$ does not cut out a monogon or a bigon.
\end{enumerate}

Curves that connect two marked points and lie entirely on the boundary of $S$ without passing through a third marked point are boundary segments. Note that boundary segments are not arcs. For any two arcs $\gamma$, $\gamma'$ in $S$, let $e(\gamma,\gamma')$  be the minimal number of crossings of arcs $\alpha$ and $\alpha'$, where $\alpha$ and $\alpha'$ range over all arcs isotopic to $\gamma$ and $\gamma'$, respectively. We say that arcs $\gamma$ and $\gamma'$ are \textit{compatible} if $e(\gamma,\gamma')=0$.\par\bigskip

An \textit{ideal triangulation} is a maximal collection of pairwise compatible arcs (together with all boundary segments). Triangulations are \textit{connected} to each other by sequences of flips. Each flip replaces a single arc $\gamma$ in a triangulation $T$ by a (unique) arc $\gamma'\neq \gamma$ that, together with the remaining arcs in $T$, forms a new triangulation.\par\bigskip

According to Schiffler and Canakci \cite{Canakci3} , Fomin, Shapiro and Thurston \cite{Fomin} associated a cluster algebra $\mathscr{A}(S,M)$ to any bordered surface with marked points $(S, M)$, and the cluster variables of $\mathscr{A}(S, M)$ are in bijection with the (tagged) arcs of $(S,M)$.\par\bigskip

The following theorem regarding relationships between cluster algebras and surface triangulations was obtained Fomin, Shapiro, and Thurston \cite{Fomin3,Fomin4}

\begin{teor} \label{teor1.6.1}  \cite{Musiker} \textit{Fix a bordered surface $(S,M)$ and let $\mathcal{A}$ be the cluster algebra associated to the signed adjacency matrix of a tagged triangulation. Then the (unlabeled) seed $\Sigma_{T}$ of $\mathcal{A}$ are in bijection with tagged triangulations $T$ of $(S,M)$, and the cluster variables are in bijection with the tagged arcs of $(S,M)$ (so we can denote each by $x_{\gamma}$, where $\gamma$ is a tagged arc).  Moreover, each seed in $\mathcal{A}$ is uniquely determined by its cluster. Furthermore, if a tagged triangulation $T'$ is obtained from another tagged triangulation $T$ by flipping a tagged arc $\gamma \in T$ and obtaining $\gamma'$, then $\Sigma_{T'}$ is obtained from $\Sigma_{T}$ by the seed mutation replacing $x_{\gamma}$ by $x_{\gamma'}$}.

\end{teor}

\subsection{Snake Graphs and Cluster Variables}

In this section, we recall the definition of a snake graph, the number of perfect matchings associated to these graphs, and the way that these concepts can be used to find out a formula for the cluster variables of a cluster algebra associated to a surface \cite{Canakci, Canakci3, Musiker}.\par\bigskip

A \textit{tile} $G$ is a square of fixed side-length in the plane whose sides are parallel or orthogonal to the fixed basis.

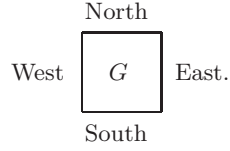
\begin{center}
\begin{picture}(60,40) 
\put(30,30){\line (1,0){30}}
\put(30,0){\line (0,1){30}}
\put(30,0){\line (1,0){30}}
\put(60,0){\line (0,1){30}}
\put(40,12){$G$}
\put(3,12){West}
\put(65,12){East.}
\put(31,35){North}
\put(31,-11){South}
\end{picture}
\par\bigskip
\end{center}

We consider a tile $G$ as a graph with four vertices and four edges in the obvious way. \par\bigskip

A \textit{snake graph} $\mathcal{G}$ is a connected graph consisting of a finite sequence of tiles $G_{1}, \dots, G_{d}$ with $d\geq 1$, such that for each $i=1, \dots, d-1$;

\begin{enumerate}
\item [(i)] $G_{i}$ and $G_{i+1}$ share exactly one edge $e_{i}$ and this edge is either the north edge of $G_{i}$ and the south edge of $G_{i+1}$ or the east edge of $G_{i}$ and the west edge of $G_{i+1}$.
\item [(ii)] $G_{i}$ and $G_{j}$ have on edge in common whenever $|i-j|\geq 2$.
\item [(iii)] $G_{i}$ and $G_{j}$ are disjoint whenever $|i-j|\geq 3$. 
\end{enumerate}

For notation, $\mathcal{G}[i, i+t]=(G_{i}, \dots, G_{i+t})$ is the subgraph of $\mathcal{G}=(G_{1}, \dots, G_{n})$, the $d-1$  edges $e_{1}, \dots, e_{d-1}$ which are contained in two tiles are called \textit{interior edges} of $\mathcal{G}$ and the other edges are  called \textit{boundary edges}. A \textit{perfect matching} $P$ of a graph $G$ is a subset of the set of edges of $G$ such that each vertex of $G$ is incident to exactly one edge in $P$. We let Match$(G)$ denote the set of all perfect matchings of the graph $G$.\par\bigskip

Let $T$ be a triangulation of a surface $(S,M)$ and let $\gamma$ be an arc in $(S,M)$ which is not in $T$. Choose an orientation on $\gamma$, let $s \in M$ be its starting point, and let $t \in M$ be its endpoint. Denote by $s=p_{0},p_{1},\dots, p_{d+1}=t$ the ordered points of intersection of $\gamma$ and $T$. For $j=1,2, \dots, d$, let $\tau_{i_{j}}$ be the arc of $T$ containing $p_{j}$, and let $\Delta_{j-1}$ and $\Delta_{j}$ be the two triangles in $T$ on either side of $\tau_{i_{j}}$. Then, for $j=1, \dots, d-1$, arcs $\tau_{i_{j}}$ and $\tau_{i_{j+1}}$ form two sides of the triangle $\Delta_{j}$ in $T$ and we define $e_{j}$ to be the third arc in this triangle.\par\bigskip

Let $G_{j}$ be the quadrilateral in $T$ that contains $\tau_{i_{j}}$ as a diagonal (a tile)  whose edges are arcs in $T$, thus, they are labeled edges. Define a sign function $f$ of the edges $e_{1}, \dots, e_{d}$ by

\begin{equation}
f(e_{j})=\begin{cases}
 +1,& \mbox{ if } e_{j} \text{ lies on the right of } \gamma \text{ when passing through } \Delta_{j},\\
 -1,& \mbox{ otherwise. }
\end{cases}
\end{equation}

The labeled snake graph $\mathcal{G}_{\gamma}=(G_{1}, \dots, G_{d})$ with tiles $G_{i}$ and sign function $f$ is called the snake graph associated to the arc $\gamma$. Each edge $e$ of $\mathcal{G}_{\gamma}$ is labeled by an arc $\tau(e)$ of the triangulation $T$. Such an arc defines the weight $x(e)$ of the edge $e$ as the cluster variable associated to the arc $\tau(e)$. Thus $x(e)=x_{\tau(e)}$.\par\bigskip

In \cite{Musiker} Musiker, Schiffler, and Williams  showed a  combinatorial formula  for cluster variables of a cluster algebra of surface type $\mathcal{A}(S,M)$ with principal coefficients $\Sigma_{T}=(\text{\textbf{x}}_{T},\text{\textbf{y}}_{T}, B_{T})$. In such a case, if $\gamma$ is a generalized arc in a triangulation  $T$ which has no self-folded triangles, and $\mathcal{G}_{\gamma}$ is its snake graph. Then  the corresponding cluster variable $x_{\gamma}$ is given by the identity 

\begin{equation}\label{equation1.9.2}
x_{\gamma}=\cfrac{1}{\text{cross}(\gamma, T)}\sum_{P \in \text{Match}(\mathcal{G}_{\gamma})}x(P).
\end{equation}

where  the sum runs over all perfect matchings of $\mathcal{G}_{\gamma}$, the summand $x(P)= \prod_{e \in P} x(e)$ is the weight of the perfect matching $P$, and $\text{cross}(T,\gamma)=\prod_{j=1}^{d}x_{\tau_{i_{j}}}$ is the product of all initial cluster variables whose arcs cross $\gamma$.\par\bigskip

A relationship between cluster variables and continued fractions is described by Schiffler and Canakci in \cite{Canakci3}, who claimed that, the numerator of a continued fraction is equal to the number of perfect matchings of the corresponding abstract snake graph, and that it can therefore be interpreted as the number of terms in the numerator of the Laurent expansion of an associated cluster variable. Thus, the Laurent polynomials of the cluster variable can be recovered from the continued fraction.

\subsection{Category of Diagonals}

In 2006 \cite{Caldero}, Caldero, Chapoton, and  Schiffler introduced the category of diagonals of a polygon with $n+3$ sides associated to a triangulation $T$, in this case, the diagonals are called \textit{roots} which can be classified as negative or positive, negative roots are those roots belonging to the triangulation $T$.\par\bigskip

The  combinatorial $\mathbb{C}$-linear additive category $C_{T}$ is  described as follows. The objects are positive integral linear combinations of positive roots, and the space of morphisms from a positive root $\alpha$ to a positive root $\alpha'$ is a quotient of the vector space over $\mathbb{C}$ spanned by pivoting paths from $\alpha$ to $\alpha'$. The subspace which defines the quotient is spanned by the so-called \textit{mesh relations}. For any couple $\alpha, \alpha'$ of positive roots such that $\alpha$ is related to $\alpha'$ by two consecutive pivoting elementary moves with distinct pivots,  the mesh relations are given by the identity $P_{v'_{2}}P_{v_{1}}=P_{v'_{1}}P_{v_{2}}$, where  $v_{1}, v_{2}$ (resp. $v'_{1} v'_{2}$) are the vertices of $\alpha$ (resp. $\alpha'$) such that $P_{v'_{1}}P_{v_{2}}=\alpha'$.\par\bigskip  

Let $T$ be a triangulation, then one can define a planar tree $t_{T}$ as follows. Its vertices are the triangles of $T$ and the edges connect adjacent triangles.  In the same way, we can define a  graph $Q_{T}$ whose vertices are the inner edges of $T$ and are related to each other by an edge, if they bound the same triangle. An orientation can be defined by using graph $Q_{T}$, in such a way that  a vertex $i$ connects a vertex $j$ (denoted $i\rightarrow j$),  if $- \alpha_{j}$ can be obtained from the diagonal $-\alpha_{i}$ by rotating anticlockwise about their common vertex. \par\bigskip

The triangulation $T$ defines a $\mathbb{C}-$linear abelian category Mod $Q_{T}$, that is, the category of modules over the quiver $Q_{T}$, such that in any triangle, the composition of two successive maps is zero. These relations are named \textit{triangle relations}.\par\bigskip

The following result regarding the category of diagonals was given by Caldero, Chapoton, and  Schiffler in \cite{Caldero}.

\begin{teor} \label{teor1.7.1}

\textit{If $T$ is a triangulation of a polygon with $n+3$ sides then there is a categorical equivalence between the category of diagonals $C_{T}$ and the category of modules over the quiver $Q_{T}$.}

\end{teor}

\section{Dyck Paths Category}

In this section, we introduce the category of Dyck paths of length $2n$. 

\subsection{Elementary Shifts}

A Dyck path is a lattice path in $\mathbb{Z}^2$  from $(0,0)$ to $(n, n)$ with steps $(1,0)$ and $(0,1)$ such that the path never passes below the line $y=x$. The number of Dyck paths of length $2n$ is equal to the nth Catalan number  \cite{Stanley}. Henceforth, Dyck words as defined in the following Remark \ref{Dyck words}  are used to denote Dyck paths.

\addtocounter{Nota}{2}

\begin{Nota}\label{Dyck words}

The set of \textit{Dyck words} is the set of words $w  \in X^{\ast}=\lbrace U,D\rbrace^{\ast}$ characterized by the following two conditions \cite{Barcucci}:

\begin{itemize}
\item for any left factor $u$ of $w$, $|u|_{U}\geq |u|_{D}$,
\item $|w|_{U}=|w|_{D}$.
\end{itemize}

where $|w|_{a}$ is the number of occurrences of the letter $a \in X$ in the word $w$ and the word $u$ is a \textit{left factor} of the word $w=uv$.\\

\end{Nota}

Let $\mathfrak{D}_{2n}$ be the set of all Dyck paths of length $2n$, let  $UWD=Uw_{1}\dots w_{n-1}D$ be a Dyck path in $\mathfrak{D}_{2n}$ with $A=\lbrace UD,DU, UU, DD\rbrace$  the set of choices in $W$.\\

The \textit{support} of $UWD$ (denoted by $\text{Supp } UWD\subseteq\{1,2,\dots, n-1\}=\textbf{n-1}$) is a set of indices such that \[\text{Supp}\hspace{0.1cm}UWD= \{ q\in\textbf{n-1} \text{ } | \text{ } w_{q}=UD \text{ or }w_{q}=UU  \text {, } 1 \leq q \leq n-1 \}. \]

A map $f:A\longrightarrow A$ such that for any $w\in A$, it holds that $f(w)= f(ab)=w^{-1}=ba$, $a, b\in\{U, D\}$ is said to be a \textit{shift}. An \textit{unitary shift} is a map $f_{i}:\mathfrak{D}_{2n}\longrightarrow \mathfrak{D}_{2n}$ such that \[f_{i}(U w_{1}\ldots w_{i-1}w_{i}w_{i+1}\ldots w_{n-1}D)=U w_{1}\ldots w_{i-1}f(w_{i})w_{i+1}\ldots w_{n-1}D.\] We will denote  a unitary shift  by a vector of maps from $\mathfrak{D}_{2n}$ to itself of the form $(1_{1}, \ldots ,1_{i-1}, f_{i}, 1_{i+1}, \ldots, 1_{n-1})$, where $1_{k}$ is the identity map associated to the $i$-th coordinate.\\

An \textit{elementary shift is a composition of unitary shifts}. A \textit{shift path} of length $m$ $UWD \longrightarrow UW_1D\longrightarrow \cdots \longrightarrow UW_{m}D \longrightarrow UVD$ from $UWD$ to $UVD$ is a composition of elementary shifts. The set of all Dyck paths in a shift path between $UWD$ and $UVD$ will be denoted by $J$. For notation, we introduce  the \textit{identity shift} as the elementary shift $(1_{1},\dots , 1_{n-1})$.\\

\textbf{Irreversibility condition}. Suppose that a map $R: \mathfrak{D}_{2n}\rightarrow \mathfrak{D}_{2n}$ is defined by the application of successive elementary shifts to a given Dyck path. Then $R$ is said to be an  \textit{irreversible relation} over $\mathfrak{D}_{2n}$   if and only if  elementary shifts transforming Dyck paths (from one to the other) are not reversible. In other words, if an elementary shift $F=f_{p_{1}}\circ \dots \circ f_{p_{q}}$ transforms a Dyck path $UWD$ into a Dyck path $UVD$ then there is not an elementary shift $F=f_{p_{1}}\circ \dots \circ f_{p_{q}}$ transforming $UVD$ into  $UWD$, for some $p,q \in \mathbb{Z}^{+}$.\\

\textbf{Shift Relation}. If there exist two paths $G \circ F$ and $G' \circ F'$ of irreversible relations (of length 2) transforming a Dyck path $UWD$ into the Dyck path $UVD$ over $R$ in the following form:

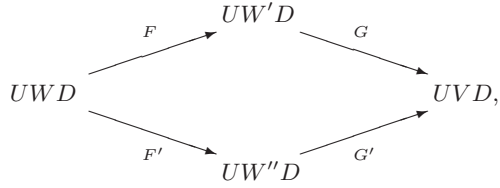
\begin{center}
\begin{picture}(170,70)
\put(0,30){$UWD$}
\put(80,60){$UW'D$}
\put(80,0){$UW''D$}
\put(160,30){$UVD$,}
\put(30,40){\vector(3,1){48}}
\put(30,26){\vector(3,-1){48}}
\put(110,10){\vector(3,1){48}}
\put(110,56){\vector(3,-1){48}}
\put(50,55){$_{F}$}
\put(50,9){$_{F'}$}
\put(130,55){$_{G}$}
\put(130,9){$_{G'}$}
\end{picture}
\end{center}

 with $W^{'}\neq W^{''}$. Then $G \circ F$  is said to be related with $G' \circ F'$ (Denoted $G \circ F \sim_{R} G' \circ F'$) whenever $G'=F$ and $G=F'$.\par\bigskip

\textbf{Category of Dyck paths of length $2n$}. As for the case of diagonals \cite{Caldero}, we can also define a  $k$-linear  additive category $(\mathfrak{D}_{2n},R)$ based on Dyck paths, in this case,  \textit{objects} are $k$-linear combinations of Dyck paths in $\mathfrak{D}_{2n}$ with  \textit{space of morphisms} from a  Dyck path $UWD$ to a Dyck path $UVD$ over $R$ being the set \[\mathrm{Hom}_{(\mathfrak{D}_{2n},R)} (UWD,UVD)= \langle  \lbrace g \text{ } \vert \text{ }  g\hspace{0.1cm} \text{is a shift path over } R \rbrace \rangle / \langle  \sim_{R} \rangle.\] The \textit{vector space} $\mathrm{Hom}_{(\mathfrak{D}_{2n},R)} (UWD,UVD)\neq0$ if and only if  there are shift paths from $UWD$ to $UVD$ and 

\begin{equation} \label{equation3.1}
\bigcap_{i \in J } \text{Supp }UW^{i}D \neq \varnothing,
\end{equation} 

for each shift path, with $UWD$ and $UVD$  in $\mathfrak{D}_{2n}$. \par\bigskip

\begin{figure}[h!]
\begin{center}
\begin{picture}(170,80)
\put(0,0){\line(1,1){30}}
\put(0,0){\line(0,1){10}}
\put(0,10){\line(1,0){10}}
\put(10,10){\line(0,1){10}}
\put(10,20){\line(1,0){10}}
\put(20,20){\line(0,1){10}}
\put(20,30){\line(1,0){10}}

\put(0,50){\line(1,1){30}}
\put(0,50){\line(0,1){10}}
\put(0,60){\line(1,0){10}}
\put(10,60){\line(0,1){20}}
\put(10,80){\line(1,0){20}}
\multiput(10,70)(3.6,0){3} {\line(1,0){2.6}} 
\multiput(20,70)(0,3.6){3} {\line(0,1){2.6}} 

\put(60,0){\line(1,1){30}}
\put(60,0){\line(0,1){20}}
\put(60,20){\line(1,0){20}}
\put(80,20){\line(0,1){10}}
\put(80,30){\line(1,0){10}}
\multiput(60,10)(3.6,0){3} {\line(1,0){2.6}} 
\multiput(70,10)(0,3.6){3} {\line(0,1){2.6}} 

\put(60,50){\line(1,1){30}}
\put(60,50){\line(0,1){20}}
\put(60,70){\line(1,0){10}}
\put(70,70){\line(0,1){10}}
\put(70,80){\line(1,0){20}}
\multiput(60,60)(3.6,0){3} {\line(1,0){2.6}} 
\multiput(70,70)(3.6,0){3} {\line(1,0){2.6}} 
\multiput(70,60)(0,3.6){3} {\line(0,1){2.6}} 
\multiput(80,70)(0,3.6){3} {\line(0,1){2.6}} 

\put(120,50){\line(1,1){30}}
\put(120,50){\line(0,1){30}}
\put(120,80){\line(1,0){30}}
\multiput(120,60)(3.6,0){3} {\line(1,0){2.6}} 
\multiput(120,70)(3.6,0){6} {\line(1,0){2.6}} 
\multiput(130,60)(0,3.6){6} {\line(0,1){2.6}} 
\multiput(140,70)(0,3.6){3} {\line(0,1){2.6}} 

\put(15,30){\vector(0,1){20}}
\put(70,30){\vector(0,1){20}}
\put(55,15){\vector(-1,0){25}}
\put(55,65){\vector(-1,0){25}}

\end{picture}
\end{center}

\caption{ Elementary shifts in $(\mathfrak{D}_{6},R)$.}\label{figure3.1}

\end{figure}
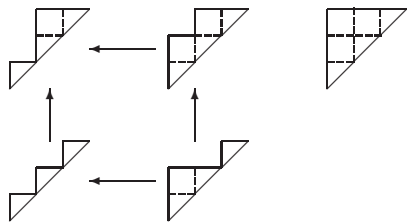

Figure \ref{figure3.1} shows the elementary shifts over $(\mathfrak{D_{6}},R)$ associated to an irreversible relation $R$ defined over the set of all Dyck paths of length $6$. And such that,

\begin{equation}
R(UWD)=\begin{cases}
  f_{1}(UWD), & \mbox{ if }  w_{1}=UD,\\
  f_{2}(UWD),& \mbox{ if } w_{2}=UD.
\end{cases}
\end{equation}

\subsection{Relations of Type $R_{j_{1}\ldots j_{m}}^{i_{1}\ldots i_{k}}$}

If $\textbf{n}=\lbrace 1,2,\ldots, n\rbrace$ is an $n$-point chain then $\mathscr{C}_{\textbf{(1,n)}}$ stands for all \textit{admissible subchains} $\mathscr{C}$ of  $\textbf{n}$  with $\text{min }\mathscr{C}=1$ and $\text{max }\mathscr{C}=n$. For instance,  $\lbrace 1,6,8 \rbrace$ and $\lbrace 1,3,8 \rbrace $ are three-point subchains contained in $\mathscr{C}_{\textbf{(1,8)}}$. The admissible subchain $\mathscr{C}=\{j_{1},\dots, j_{m},i_{1},\dots, i_{k}\}\subseteq\textbf{n}$ must satisfy the following constraints for $1\leq r\leqslant k$ and $1\leq s \leqslant m$:

\begin{itemize}
 \item If $i_{1}=1$ and $k=m$ then $i_{1}<j_{1}< \dots<i_{k}<j_{m}=n$.\\
\item If $i_{1}=1$ and $k=m+1$ then $i_{1}<j_{1}< \dots<i_{k}<j_{m}<i_{k}=n$.\\
\item If  $j_{1}=1$ and $k=m$ then $j_{1}<i_{1}< \dots < j_{m}<i_{k}=n$.\\
\item If $j_{1}=1$ and $m=k+1$ then $j_{1}<i_{1}< \dots < j_{m}<i_{k}<j_{m}=n$.
\end{itemize}

Let $\sigma:\{i_{1},j_{1}\} \rightarrow \{0,1\}$  be a map such that $\sigma(i_{1})=1$ and $\sigma(j_{1})=0$. For $a \in \{ i_{1},j_{1}\}$, we assume  $i_{r} (i_{r+\sigma(a)}) \in \{i_{1},\dots, i_{k}\}$ and $ j_{r+1-\sigma(a)} (j_{r}) \in \{j_{1},\dots , j_{m}\}$. The orientation between the interval $\lbrack i_{r},j_{r+1-\sigma(a)}\rbrack $ ($\lbrack j_{r},i_{r+\sigma(a)}\rbrack$) over $\mathbb{Z}^{+}$ is from left to right, denoted $\overrightarrow{[a, b]}$ (right to left, denoted  $\overleftarrow{[a ,b]}$). The following words are defined by using intervals:

\begin{itemize}

\item $w_{t}= \min \lbrace  \text{ } w_{s} \text{ }\vert \text{ } w_{i_{r}} \leq w_{s} \leq w_{j_{r+1-\sigma(a)}}, w_{s}=UD \text{ } \rbrace$   $\big (  w_{t}= \max \lbrace  \text{ } w_{s} \text{ }\vert \text{ } w_{j_{r}} \leq w_{s} \leq w_{i_{r+\sigma(a)}}, w_{s}=UD \text{ } \rbrace  \big )$,

\item $w_{p}= \min \lbrace  \text{ } w_{s} \text{ }\vert \text{ } w_{t} < w_{s} \leq w_{j_{r+1 -\sigma(a)}}, w_{s}=DU \text{ } \rbrace$  $\big (  w_{p}= \max \lbrace  \text{ } w_{s} \text{ }\vert \text{ } w_{j_{s}} \leq w_{s} < w_{t}, w_{s}=DU \text{ } \rbrace  \big )$.

\end{itemize}  

We introduce the following elementary shifts:

\begin{enumerate}[ES1.]

\item If $w_{s}=UD$ for all $w_{s} \in \lbrack i_{r},j_{r+1-\sigma(a)}\rbrack$ ($\lbrack j_{r},i_{r+\sigma(a)}\rbrack$), \[ \overrightarrow{\lbrack j_{r-\sigma(a)}, i_{r} \rbrack} \overleftarrow{\lbrack i_{r},j_{r+1-\sigma(a)}\rbrack}  \overrightarrow{\lbrack j_{r+1-\sigma(a)}, i_{r+1} \rbrack}, \] \[ ( \overleftarrow{\lbrack i_{r+\sigma(a)-1},j_{r}\rbrack} \overrightarrow{\lbrack j_{r}, i_{r+\sigma(a)} \rbrack} \overleftarrow{\lbrack i_{r+\sigma(a)},j_{r+1} \rbrack}),\] then 
 \[g(UWD)= f_{j_{r+1-\sigma(a)}} \circ \dots  \circ f_{i_{r}}(UWD),\] if there exists $s \in \mathbb{Z}^{+}$ such that $j_{r-\sigma(a)}\leq s \leq i_{r}$, $| s-j_{r}| > 1$, $w_{x}=UD$ if $s \leq x \leq i_{r}$ over $\lbrack j_{r-\sigma(a)}, i_{r} \rbrack$ and 

\begin{equation}
 w_{y}=\begin{cases}
 UD, & \mbox{if $y=j_{r+1-\sigma(a)}$,}\\
 DU,& \mbox{otherwise,}
\end{cases} 
\end{equation}

over $[ j_{r+1-\sigma(a)}, i_{r+1} ]$ for $j_{r+1-\sigma(a)}\neq n-1$ or the first condition over $[ j_{r-\sigma(a)}, i_{r} ]$ for $j_{r+1-\sigma(a)}=n-1$.

  \[\Bigg( g(UWD)=f_{i_{r+\sigma(a)}} \circ \dots \circ f_{j_{r}} (UWD),\] if there exists $s \in \mathbb{Z}^{+}$ such that $i_{r+\sigma(a)} \leq s \leq j_{r+1}$, $|s-i_{r+\sigma(a)}|>1$, $w_{x}=UD$ if $i_{r+\sigma(a)} \leq x \leq s$ over $\lbrack i_{r+\sigma(a)},j_{r+1} \rbrack$ and 

 \begin{equation}
 w_{y}=\begin{cases}
   UD, & \mbox{if $y=j_{r}$,}\\
 DU,& \mbox{otherwise,}
\end{cases}
\end{equation}  

 over $\lbrack i_{r+\sigma(a)-1},j_{r}\rbrack$ for $j_{r}\neq 1$ or the first condition over $\lbrack i_{r+\sigma(a)},j_{r+1} \rbrack$ for $j_{r}=1$\Bigg), with $i_{r}\neq 1$ ($i_{r+\sigma(a)}\neq n-1$).

 \item  If $t=1$ or $n-1$ then $g(UWD)=f_{t}(UWD)$. Then there are not elementary shifts associated to elements of a subchain different from $1$ and $n-1$.

\item If  $w_{i_{r}}< w_{t} < w_{j_{r +1 -\sigma(a)}}$ ($w_{j_{r}}<w_{t}<w_{i_{r+\sigma(a)}}$) then $g(UWD)=f_{t}(UWD)$.

\item If $w_{p}=w_{j_{r+1-\sigma(a)}} (w_{j_{r}})$  then 

\[ g(UWD)=\begin{cases}
   f_{i_{r+1}}\circ \dots  \circ f_{j_{r+1-\sigma(a)}}(UWD) & \mbox{if $j_{r+1-\sigma(a)} \neq n-1  $,}\\
 f_{j_{r+1-\sigma(a)}}(UWD) & \mbox{if $j_{r+1-\sigma(a)} = n-1$.}
\end{cases}\] 

\[ \Bigg (g(UWD)=\begin{cases}
   f_{i_{r+\sigma(a)-1}}\circ \dots  \circ f_{j_{r}}(UWD) & \mbox{if $j_{r} \neq 1  $,}\\
 f_{j_{r}}(UWD) & \mbox{if $j_{r} = 1.$}
\end{cases} \Bigg)\]

\item If $w_{t}< w_{p}< w_{j_{r+1-\sigma(a)}}$ ($w_{j_{r}}<w_{p}<w_{t}$) then $g(UWD)=f_{p}(UWD)$.

\end{enumerate}

\par\bigskip

For a given subchain  $\mathscr{C}= \{j_{1},\dots, j_{m},i_{1},\dots, i_{k}\} \subseteq\textbf{n-1}$, two Dyck paths $D$ and $D'$ of length $2n$ are said to be \textit{related by a relation of type} $R_{j_{1}\ldots j_{m}}^{i_{1}\ldots i_{k}}$ if there is an elementary shift $\mathrm{ES}i$, $1\leq i\leq 5$ which transforms either $D$ into $D'$ or $D'$ into $D$.

\addtocounter{prop}{3}

 \begin{prop}\label{prop3.1}

\textit{The relation $R_{j_{1}\ldots j_{m}}^{i_{1}\ldots i_{k}}$ is irreversible.}

 \end{prop}

 \textbf{Proof.} Suppose that there is an elementary shift  $f_{r_{1}}\circ \dots \circ f_{r_{t}}$ from a Dyck path $UWD$ to a Dyck path $UVD$ and that there is an elementary shift $f_{r_{1}}\circ \dots \circ f_{r_{t}}$ from $UVD$ to a $UWD$, then we have five cases: 

 \begin{enumerate}[(i)]

     \item  If $f_{r_{1}}\circ \dots \circ f_{r_{t}}$ arises  from ES1 over $[i_{r},j_{r+1-\sigma(a)}]$.  Shifts ES2, ES3 and ES5 allow to conclude that from $UVD$ to a $UWD$, $f_{t}=f_{j_{r+1-\sigma(a)}}\circ \dots \circ f_{i_{r}}$ or $f_{p}=f_{j_{r+1-\sigma(a)}}\circ \dots \circ f_{i_{r}}$ and this is a contradiction. If ES1 is an elementary shift from $UVD$ to  $UWD$, then two cases arise: If $j_{r+1-\sigma(a)} \neq n-1$, thus $UVD$ equals

  { \scriptsize  \[ Uv_{1}\dots v_{j_\sigma(a)} \dots \underbrace{v_{j_{s}} \dots  v_{i_{r}-1}}_{UD} \underbrace{v_{i_{r}} \dots v_{j_{r+1-\sigma(a)}} }_{UD} \underbrace{(v_{j_{r+1-\sigma(a)}+1}) \dots v_{i_{r+1}} }_{DU} v_{r_{i+1}+1}\dots v_{n-1}D, \]}

     it turns out that $f_{j_{r+1-\sigma(a)}}\circ \dots \circ f_{i_{r}}(UVD)$ has the form

     {\scriptsize  \[ Uw_{1}\dots w_{j_\sigma(a)} \dots \underbrace{w_{j_{s}} \dots  (w_{i_{r}}-1)}_{UD} \underbrace{w_{i_{r}} \dots w_{j_{r+1-\sigma(a)}} }_{DU} \underbrace{(w_{j_{r+1-\sigma(a)}+1}) \dots w_{i_{r+1}} }_{DU}w_{r_{i+1}+1}\dots v_{n-1}D, \]} which is a contradiction. If $j_{r+1-\sigma(a)}=n-1$, $UVD$ is equal to

 { \scriptsize  \[ Uv_{1}\dots v_{j_\sigma(a)} \dots \underbrace{v_{j_{s}} \dots  v_{i_{r}-1}}_{UD} \underbrace{v_{r} \dots v_{j_{r+1-\sigma(a)}} }_{UD}D, \]} and  $f_{j_{r+1-\sigma(a)}}\circ \dots \circ f_{i_{r}}(UVD)$ has the shape { \scriptsize \[ Uw_{1}\dots w_{j_\sigma(a)} \dots \underbrace{w_{j_{s}} \dots  w_{i_{r}-1}}_{UD} \underbrace{w_{r} \dots w_{j_{r+1-\sigma(a)}} }_{DU}D,\]} again a contradiction. We also get a contradiction if a elementary shift is done by using ES4 from $UVD$ to a $UWD$, indeed, in these cases it holds that, if $j_{r+1-\sigma(a)}= n-1$, there are $t$ and $p$ such that  $ p=j_{r+1-\sigma(a)} <t  \leq i_{r+1}$ and $UVD$ is equal to

 { \scriptsize  \[ Uv_{1}\dots  v_{i_{r}-1} \underbrace{v_{i_{r}} \dots v_{j_{r+1-\sigma(a)}} }_{DU} \underbrace{v_{j_{r+1-\sigma(a)+1}} \dots v_{t} }_{UD}v_{t+1}\dots v_{n-1}D, \]} and $f_{j_{r+1-\sigma(a)}}\circ \dots \circ f_{i_{r}}(UVD)$ is { \scriptsize \[ Uw_{1}\dots  w_{i_{r}-1}  \underbrace{w_{i_{r}} \dots w_{j_{r+1-\sigma(a)}} }_{UD} \underbrace{w_{j_{r+1-\sigma(a)+1}} \dots w_{t} }_{UD} w_{t+1 }\dots w_{n-1}D. \]} If $v_{j_{r+1-\sigma(a)}}=n-1$ $f_{r+1-\sigma(a)}= f_{r+1- \sigma(a)}\circ \dots \circ f_{i_{r}}$  but this is a contradiction.

 \item If $f_{r_{1}}\circ \dots \circ f_{r_{t}}$ arises from ES2 over $[i_{1},j_{1}]$ then we cannot use elementary shifts defined in cases ES1, ES4, ES5 or ES3, provided that, $i_{1}\neq1$, $t\neq p$ or $1<t<j_{1}$. Therefore, ES2 guarantees the existence of a walk from $UVD$ to $UWD$ such that; \[U \underbrace{v_{1}}_{UD} \dots v_{j_{1}} \dots v_{n-1}D,\] and $f_{1}(UWD)$ has the form \[U \underbrace{w_{1}}_{DU} \dots w_{j_{1}} \dots w_{n-1}D,\] which is a contradiction (if $t=n-1$, the proof is dual).

 \item If $f_{r_{1}}\circ \dots \circ f_{r_{t}}$ arises  from ES3 over $[i_{r},j_{r+1-\sigma(a)}]$, provided that,  $i_{r}<t<p<j_{r+1-\sigma(a)}$, we conclude that it is not possible to use ES1, ES2, ES4 nor ES5. In the case of ES3 from $UVD$ to a $UWD$, $UVD$ equals \[Uv_{1}\dots \underbrace{v_{i_{r}} \dots v_{t-1}}_{DU} \underbrace{v_{t}}_{UD}\dots  v_{j_{r+1-\sigma(a)}} \dots v_{n-1}D,\] and $f_{t}(UVD)$ has the shape \[Uw_{1}\dots \underbrace{w_{i_{r}} \dots w_{t}}_{DU}w_{t+1}\dots  w_{j_{r+1-\sigma(a)}} \dots w_{n-1}D,\] but this is a contradiction.

 \item If $f_{r_{1}}\circ \dots \circ f_{r_{t}}$ arises  from ES4 over $[i_{r},j_{r+1-\sigma(a)}]$, provided that $t<p$, we do not use ES2, ES3 nor ES5. If  $j+1-\sigma(a) = n-1$, we cannot use ES1. If $j+1-\sigma(a) \neq n-1$ we can use ES1 from $UVD$ to a $UWD$ (Note that, it is not necessary with $v_{m}=UD$ for all $s\in [j_{r+1-\sigma(a)}+1,i_{r+1}]$) $UVD$ is equal to { \scriptsize \[Uv_{1}\dots v_{i_{r}-1}\underbrace{v_{i_{r}} \dots v_{t} \dots v_{p-1}}_{DU} \underbrace{v_{p}v_{j_{r+1-\sigma(a)}+1} \dots v_{i_{r+1}} }_{UD} v_{i_{r+1}+1}\dots v_{n-1}D,\] } it turns out that $g(UVD)$ has the form { \scriptsize \[Uw_{1}\dots w_{i_{r}-1}\underbrace{w_{i_{r}} \dots w_{t} \dots w_{p-1} w_{p}w_{j_{r+1-\sigma(a)}+1} \dots w_{i_{r+1}} }_{DU} w_{i_{r+1}+1}\dots w_{n-1}D,\]}which is a contradiction. Using ES5 from $UVD$  to  $UWD$, if $j_{r+1-\sigma(a)} \neq n-1$, $UVD$ is equal to  \[Uv_{1} \dots v_{i_{r}} \dots \underbrace{v_{t}\dots v_{p-1}}_{UD}\underbrace{v_{p}}_{DU}v_{j_{r+1-\sigma(a)}+1} \dots v_{i_{r+1}}v_{i_{r+1}+1}\dots v_{n-1}D,\] and $UWD$ has the shape \[Uw_{1} \dots w_{i_{r}} \dots \underbrace{w_{t}\dots w_{p}}_{UD}\underbrace{w_{j_{r+1-\sigma(a)}+1} \dots w_{i_{r+1}}}_{f(ab)}w_{i_{r+1}+1}\dots w_{n-1}D,\] again a contradiction. If $j_{r+1-\sigma(a)} = n-1$, $UVD$ is equal to \[Uv_{1}\dots v_{i_{r}}\dots v_{t-1}\underbrace{v_{t}\dots v_{p-1}}_{UD}\underbrace{v_{p}}_{DU}D,\] it turns out that $UWD$ has the shape \[Uw_{1}\dots w_{i_{r}}\dots w_{t-1}\underbrace{w_{t}\dots v_{p}}_{UD}D,\] this is a contradiction.

 \item If $f_{r_{1}}\circ \dots \circ f_{r_{t}}$ arises  from ES5 over $[i_{r},j_{r+1-\sigma(a)}]$. Then we cannot use ES1, ES2, ES3 nor ES4, because $f_{p}\neq f_{j_{r+1-\sigma(a)}}\circ \dots \circ f_{i_{r}}$ and $t<p$. Using ES5 from $UVD$ to a $UWD$, we observe that $UVD$ is equal to \[Uv_{1}\dots v_{i_{r}} \dots v_{t-1} \underbrace{v_{t} \dots v_{p-1}}_{UD} \underbrace{v_{p}}_{DU}\dots  v_{j_{r+1-\sigma(a)}} \dots v_{n-1}D,\] and $f_{p}(UWD)$ has the form \[Uw_{1}\dots w_{i_{r}} \dots w_{t-1} \underbrace{w_{t} \dots w_{p}}_{UD} v_{p+1}\dots  v_{j_{r+1-\sigma(a)}} \dots v_{n-1}D,\]  again this is a contradiction.

\end{enumerate}

 Taking into account that if $f_{r_{1}}\circ \dots \circ f_{r_{t}}$ arises from $ES1$, $ES2$, $ES3$, $ES4$ and $ES5$ over $[i_{r},j_{r+\sigma(a)}]$ then same arguments as described above applied dually allow to conclude the proposition. We are done. $\hfill\square$
 
\subsection{$\mathbb{A}_{n-1}$-Dyck Paths Categories} 

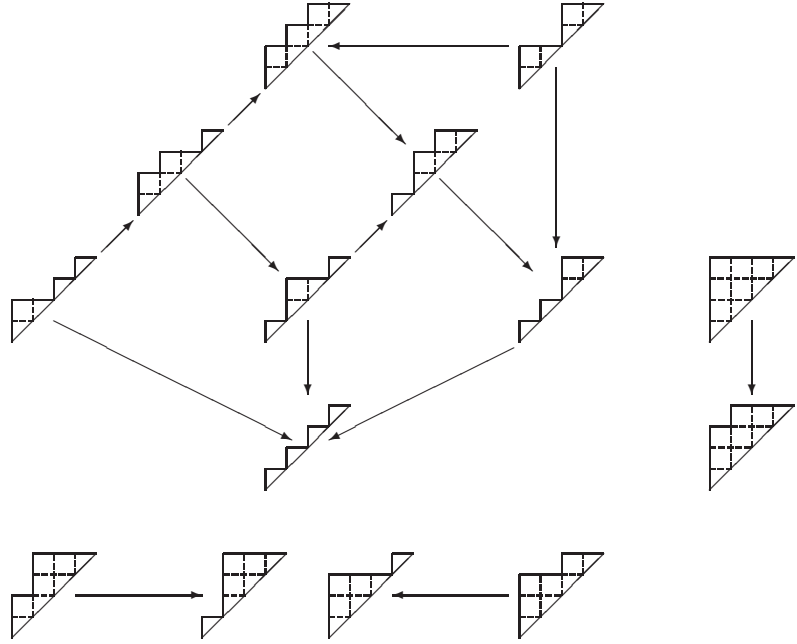
\begin{figure}[h!]

\begin{center}

\begin{picture}(288,260)

\put(0,0){\line(1,1){32}}
\put(0,0){\line(0,1){16}}
\put(0,16){\line(1,0){8}}
\put(8,16){\line(0,1){16}}
\put(8,32){\line(1,0){24}}
\multiput(0,8)(3.1,0){3} {\line(1,0){2.1}} 
\multiput(8,8)(0,3.1){3} {\line(0,1){2.1}} 
\multiput(8,16)(3.1,0){3} {\line(1,0){2.1}}
\multiput(16,16)(0,3.3){5} {\line(0,1){2.3}} 
\multiput(8,24)(3.3,0){5} {\line(1,0){2.3}}
\multiput(24,24)(0,3.1){3} {\line(0,1){2.1}} 

\put(72,0){\line(1,1){32}}
\put(72,0){\line(0,1){8}}
\put(72,8){\line(1,0){8}}
\put(80,8){\line(0,1){24}}
\put(80,32){\line(1,0){24}}
\multiput(80,16)(3.1,0){3} {\line(1,0){2.1}} 
\multiput(80,24)(3.3,0){5} {\line(1,0){2.3}} 
\multiput(88,16)(0,3.3){5} {\line(0,1){2.3}}
\multiput(96,24)(0,3.1){3} {\line(0,1){2.1}} 

\put(120,0){\line(1,1){32}}
\put(120,0){\line(0,1){24}} 
\put(120,24){\line(1,0){24}}
\put(144,24){\line(0,1){8}}
\put(144,32){\line(1,0){8}}
\multiput(120,8)(3.1,0){3} {\line(1,0){2.1}} 
\multiput(120,16)(3.3,0){5} {\line(1,0){2.3}}
\multiput(128,8)(0,3.3){5} {\line(0,1){2.3}}
\multiput(136,16)(0,3.1){3} {\line(0,1){2.1}}

\put(192,0){\line(1,1){32}}
\put(192,0){\line(0,1){24}} 
\put(192,24){\line(1,0){16}}
\put(208,24){\line(0,1){8}}
\put(208,32){\line(1,0){16}}
\multiput(192,8)(3.1,0){3} {\line(1,0){2.1}} 
\multiput(208,24)(3.1,0){3} {\line(1,0){2.1}} 
\multiput(192,16)(3.3,0){5} {\line(1,0){2.3}}
\multiput(200,8)(0,3.3){5} {\line(0,1){2.3}}
\multiput(208,16)(0,3.1){3} {\line(0,1){2.1}}
\multiput(216,24)(0,3.1){3} {\line(0,1){2.1}}

\put(96,56){\line(1,1){32}}
\put(96,56){\line(0,1){8}}
\put(96,64){\line(1,0){8}}
\put(104,64){\line(0,1){8}}
\put(104,72){\line(1,0){8}}
\put(112,72){\line(0,1){8}}
\put(112,80){\line(1,0){8}}
\put(120,80){\line(0,1){8}}
\put(120,88){\line(1,0){8}}

\put(264,56){\line(1,1){32}}
\put(264,56){\line(0,1){24}}
\put(264,80){\line(1,0){8}}
\put(272,80){\line(0,1){8}}
\put(272,88){\line(1,0){24}}
\multiput(264,64)(3.1,0){3} {\line(1,0){2.1}} 
\multiput(264,72)(3.3,0){5} {\line(1,0){2.3}}
\multiput(272,80)(3.3,0){5} {\line(1,0){2.3}}
\multiput(272,64)(0,3.3){5} {\line(0,1){2.3}}
\multiput(280,72)(0,3.3){5} {\line(0,1){2.3}}
\multiput(288,80)(0,3.1){3} {\line(0,1){2.1}}

\put(0,112){\line(1,1){32}}
\put(0,112){\line(0,1){16}}
\put(0,128){\line(1,0){16}}
\put(16,128){\line(0,1){8}}
\put(16,136){\line(1,0){8}}
\put(24,136){\line(0,1){8}}
\put(24,144){\line(1,0){8}}
\multiput(0,120)(3.1,0){3} {\line(1,0){2.1}} 
\multiput(8,120)(0,3.1){3} {\line(0,1){2.1}}

\put(96,112){\line(1,1){32}}
\put(96,112){\line(0,1){8}}
\put(96,120){\line(1,0){8}}
\put(104,120){\line(0,1){16}}
\put(104,136){\line(1,0){16}}
\put(120,136){\line(0,1){8}}
\put(120,144){\line(1,0){8}}
\multiput(104,128)(3.1,0){3} {\line(1,0){2.1}} 
\multiput(112,128)(0,3.1){3} {\line(0,1){2.1}}

\put(192,112){\line(1,1){32}}
\put(192,112){\line(0,1){8}}
\put(192,120){\line(1,0){8}}
\put(200,120){\line(0,1){8}}
\put(200,128){\line(1,0){8}}
\put(208,128){\line(0,1){16}}
\put(208,144){\line(1,0){16}}
\multiput(208,136)(3.1,0){3} {\line(1,0){2.1}} 
\multiput(216,136)(0,3.1){3} {\line(0,1){2.1}}

\put(264,112){\line(1,1){32}}
\put(264,112){\line(0,1){32}}
\put(264,144){\line(1,0){32}}
\multiput(264,120)(3.1,0){3} {\line(1,0){2.1}} 
\multiput(264,128)(3.3,0){5} {\line(1,0){2.3}}
\multiput(264,136)(3.1,0){8} {\line(1,0){2.1}} 
\multiput(272,120)(0,3.1){8} {\line(0,1){2.1}}
\multiput(280,128)(0,3.3){5} {\line(0,1){2.3}}
\multiput(288,136)(0,3.1){3} {\line(0,1){2.1}}

\put(48,160){\line(1,1){32}}
\put(48,160){\line(0,1){16}}
\put(48,176){\line(1,0){8}}
\put(56,176){\line(0,1){8}}
\put(56,184){\line(1,0){16}}
\put(72,184){\line(0,1){8}}
\put(72,192){\line(1,0){8}}
\multiput(48,168)(3.1,0){3} {\line(1,0){2.1}}
\multiput(56,176)(3.1,0){3} {\line(1,0){2.1}}
\multiput(56,168)(0,3.1){3} {\line(0,1){2.1}}
\multiput(64,176)(0,3.1){3} {\line(0,1){2.1}}

\put(144,160){\line(1,1){32}}
\put(144,160){\line(0,1){8}}
\put(144,168){\line(1,0){8}}
\put(152,168){\line(0,1){16}}
\put(152,184){\line(1,0){8}}
\put(160,184){\line(0,1){8}}
\put(160,192){\line(1,0){16}}
\multiput(152,176)(3.1,0){3} {\line(1,0){2.1}}
\multiput(160,184)(3.1,0){3} {\line(1,0){2.1}}
\multiput(160,176)(0,3.1){3} {\line(0,1){2.1}}
\multiput(168,184)(0,3.1){3} {\line(0,1){2.1}}

\put(96,208){\line(1,1){32}}
\put(96,208){\line(0,1){16}}
\put(96,224){\line(1,0){8}}
\put(104,224){\line(0,1){8}}
\put(104,232){\line(1,0){8}}
\put(112,232){\line(0,1){8}}
\put(112,240){\line(1,0){16}}
\multiput(96,216)(3.1,0){3} {\line(1,0){2.1}}
\multiput(104,224)(3.1,0){3} {\line(1,0){2.1}}
\multiput(112,232)(3.1,0){3} {\line(1,0){2.1}}
\multiput(104,216)(0,3.1){3} {\line(0,1){2.1}}
\multiput(112,224)(0,3.1){3} {\line(0,1){2.1}}
\multiput(120,232)(0,3.1){3} {\line(0,1){2.1}}

\put(192,208){\line(1,1){32}}
\put(192,208){\line(0,1){16}}
\put(192,224){\line(1,0){16}}
\put(208,224){\line(0,1){16}}
\put(208,240){\line(1,0){16}}
\multiput(192,216)(3.1,0){3} {\line(1,0){2.1}}
\multiput(208,232)(3.1,0){3} {\line(1,0){2.1}}
\multiput(200,216)(0,3.1){3} {\line(0,1){2.1}}
\multiput(216,232)(0,3.1){3} {\line(0,1){2.1}}

\put(24,16){\vector(1,0){48}}
\put(188,16){\vector(-1,0){44}}
\put(16,120){\vector(2,-1){90}}
\put(190,110){\vector(-2,-1){70}}
\put(112,120){\vector(0,-1){28}}
\multiput(34,146)(96,0){2} {\vector(1,1){12}}
\put(82,194){\vector(1,1){12}}
\multiput(66,174)(96,0){2} {\vector(1,-1){35}}
\put(114,222){\vector(1,-1){35}}
\put(188,224){\vector(-1,0){68}}
\put(206,216){\vector(0,-1){68}}
\put(280,120){\vector(0,-1){28}}

\end{picture}
\caption{Elementary shifts in an $\mathbb{A}_{3}$-Dyck paths category.}\label{figure3.2}

\end{center}

\end{figure}

For $n\geq2$ fixed, the $\mathbb{A}_{n-1}$-Dyck paths category is a category of Dyck paths $(\mathfrak{D}_{2n},R)$  where $R$ is a relation of type $R_{j_{1}\ldots j_{m}}^{i_{1}\ldots i_{k}}$ as described before. As an example we let  $(\mathfrak{D}_{8},R^{1}_{3})$ denote the $\mathbb{A}_{3}$-Dyck paths category with the admissible subchain $1<3$. Figure \ref{figure3.2} shows all the elementary shifts of $(\mathfrak{D}_{8},R^{1}_{3})$.

\par\bigskip

We let $S$ denote the set of all Dyck paths with exactly $n-1$ peaks. The following propositions and lemmas describe some properties of the set $S$.

\begin{prop} \label{prop3.2}

\textit{Let $UWD$ be a Dyck path of length $2n$, then $UWD \in S$  if and only if there is a unique sequence $w_{l}w_{l+1}\dots w_{r-1}w_{r}$ such that}

\begin{equation}\label{equation3.2}
w_{i}=\begin{cases}
  UD, & \text{if}\hspace{0.1cm} l \leq i \leq r,\\
 DU, & otherwise.
\end{cases}
\end{equation} 

\end{prop}

\textbf{Proof.} Firstly, let $\delta$ be a map $\delta: \{ ), ( \}\rightarrow \{ U,D \}$ where left bracket is associated to the letter $U$ and  right bracket is associated to the letter $D$, suppose  $UWD \in S$, then there exist bracket-subchains such that $UWD$ can be written in the following form \[\text{ }( \text{ } \underbrace{ )( }_{1} \underbrace{ )( }_{2}\dots \underbrace{ )( }_{l-2}\underbrace{ )( }_{l-1} \text{ }( \text { } \underbrace{ ( \text{ }) }_{l} \dots \underbrace{ ( \text{ }) }_{r}  \text{ }) \text{ } \underbrace{ )( }_{r+1} \dots \underbrace{ )( }_{n-2} \underbrace{ )( }_{n-1} \text{ }) \text{ },\]  therefore $w_{i}= UD$ if $l \leq i \leq r$ and $w_{i}=DU$. On the other hand, suppose $UWD$ has a unique subsequence $w_{l}w_{l+1}\dots w_{r-1}w_{r}$ that satisfies identity (\ref{equation3.2}), then if we apply $\delta^{-1}$ to $UWD$, the sequence \[ \underbrace{ ( \text{ } )}_{1} \underbrace{ ( \text{ } ) }_{2} \dots \underbrace{ ( \text{ } ) }_{l-1} \text{ } ( \text{ } \underbrace{ ( \text{ } ) }_{l} \dots \underbrace{ ( \text{ } ) }_{r} \text{ } ) \text{ } \underbrace{ ( \text{ } ) }_{r+1} \dots \underbrace{ ( \text{ } ) }_{n-2} \underbrace{ ( \text{ } ) }_{n-1},\] is obtained, therefore $UWD \in S$. We are done.  $\hfill\square$

\addtocounter{lema}{5}

\begin{lema}\label{lema3.1}

\textit{Let $UWD$ be a Dyck path in $S$, and integers $r, l$ defined as in $\mathrm{Proposition \hspace{0.1cm}\ref{prop3.2}}$ with $|r-l|>0$, then there exists an elementary shift from $UWD$ to another Dyck path with exactly $n-1$ peaks.}

\end{lema}

\textbf{Proof} Let $UWD$ be a Dyck path in $S$, let $l$ and $r$ be positive integers such that $w_{m}=UD$ for $l\leq m \leq r $. Let $l \in [i_{r},j_{r+1-\sigma(a)}]$, we have the following cases:

\begin{enumerate}[(1)]

    \item If $l=i_{r}=1$, then  \[g(UWD)= U\underbrace{f(w_{1})}_{DU} \underbrace{ w_{2} \dots w_{r}}_{UD}w_{r+1} \dots w_{n-1}D  \in S.\] 
    \item If $l=i_{r} \neq 1$, then there is $p=l-1$ over $[j_{r- \sigma(a)},i_{r}]$ such that  \[g(UWD) = Uw_{1} \dots \underbrace{f(w_{p})w_{l}  \dots w_{m}}_{UD}  \dots w_{n-1}D \in S.\]  
    \item If $i_{r}< l < j_{r+1- \sigma(a)}$, then   \[g(UWD)=Uw_{1}\dots \underbrace{f(w_{l})}_{DU}\underbrace{w_{l+1}\dots w_{r} }_{UD}w_{r+1}\dots w_{n-1}D \in S.\]  
    \item If $l=j_{r+1-\sigma(a)}$ and $|l-r|>0$, then $r \in [i_{r_{1}},j_{r_{1}+1-\sigma(a)}]$ with $|r_{1}-r|>0$ and the following cases hold:

    \begin{enumerate}[(4.1)]

        \item If $i_{r_{1}}\leq r < j_{r_{1}+1-\sigma(a)}$, there is $p=r+1$ such that, if $p\neq j_{r_{1}+1-\sigma(a)}$ then \[g(UWD)=uw_{1}\dots \underbrace{w_{l}\dots w_{r}f(w_{p})}_{UD}  \dots w_{n-1}D \in S,\]   if $p=j_{r_{1}+1-\sigma(a)}=n-1$, then \[g(UWD)=Uw_{1} \dots \underbrace{w_{l}\dots w_{r}f(w_{p})}_{UD}D \in S,\] or if $p=j_{r_{1}+1-\sigma(a)}\neq n-1$ then \[g(UWD)=Uw_{1}\dots \underbrace{w_{l} \dots w_{r}f(w_{p}) \dots f(w_{i_{r_{1}+1}})}_{UD}\dots w_{n-1}D \in S.\]
        \item If $r=j_{r_{1}+1-\sigma(a)}$ \[g(UWD)= Uw_{1}\dots \underbrace{w_{l}\dots w_{i_{r_{1}-1}}}_{UD}\underbrace{f(w_{i_{r_{1}}})\dots f(w_{r})}_{DU} \dots D \in S.\]
        \item Now, if $|r_{1}-r|>1$ or $r_{1}=r+1$ and $r>i_{r+1}+2$ then \[g(UWD)=Uw_{1}\dots \underbrace{f(w_{l})\dots f(w_{i_{r+1}})}_{DU} \underbrace{w_{i_{r+1}+1}\dots  w_{r}}_{UD}\dots D \in S.\]    

        \end{enumerate}

    For $r \in [j_{r_{1}+1-\sigma(a)},i_{r_{1}+1}]$ with $|r_{1}-r|\geq 0$ we have that:

    \begin{enumerate}[(4.4)]

        \item If $s=t=i_{r_{1}+1}=n-1$, then \[g(UWD)=Uw_{1} \dots \underbrace{w_{l} \dots w_{r-1}}_{UD}\underbrace{f(w_{r})}_{DU}D \in S. \] On the other hand, if $s=t=i_{r_{1}+1} \neq n-1$, then there is $p\in [i_{r_{1}+1}, j_{r_{1}+2-\sigma(a)}]$  satisfying first condition of (4.1). Thus, if $j_{r_{1}+1-\sigma(a)}< s < i_{r_{1}+1}$, it holds that \[g(UWD)=Uw_{1}\dots \underbrace{w_{l}\dots w_{r-1}}_{UD}\underbrace{f(w_{r})}_{DU} \dots w_{n-1}D \in S.\]

        \item [(4.5)]If $s=j_{r_{1}+1-\sigma(a)}$ then $|r_{1}-r|>0$ (If $|r_{1}-r|=0$, $|l-f|=0$ which is a contradiction) \[g(UWD)=Uw_{1}\dots \underbrace{w_{l}\dots w_{i_{r_{1}-1}}}_{UD}\underbrace{f(w_{i_{r_{1}}})\dots f(w_{r})}_{DU}w_{r+1} \dots w_{n-1}D \in S.\]

        \item [(4.6)] Now,  suppose that in $UWD$ $|r_{1}-r|>0$, then it satisfies the first condition in (4.3).

    \end{enumerate}

\end{enumerate}

In case that $l \in [j_{r},i_{r+\sigma(a)}]$, we have the following cases:

\begin{enumerate} [(5)]

    \item If $jr < l \leq i_{r}+\sigma(a)$, then there exists $p=l+1$ such that, if $p \neq j_{r}$ then \[g(UWD)= Uw_{1}\dots w_{j_{r}}\dots  \underbrace{f(w_{p})w_{l}\dots w_{r}}_{UD} \dots w_{n-1}D \in S. \] Note that, if $p=j_{r}=1$ then \[g(UWD)=U \underbrace{f(w_{p})w_{l}\dots w_{r}}_{UD} \dots w_{n-1}D \in S,\] or if $p=j_{r}\neq 1$, then \[g(UWD)=Uw_{1}\dots \underbrace{f(w_{i_{r}-1+\sigma(a)}) \dots f(w_{p})w_{l} \dots w_{r}}_{UD} \dots w_{n-1}D  \in S.\]

    \item [(6)] If $l=j_{r}$ and $|l-r|>0$, then $r \in [j_{r_{1}},i_{r_{1}+\sigma(a)}]$ with $|r_{1}-r|\geq 0$, then the following cases hold: 

    \begin{enumerate}[(6.1)]

    \item If $j_{r_{1}}+2 \leq r \leq i_{r_{1}+\sigma(a)}$, then there exists $p$ satisfying (4.4).

    \item If $j_{r_{1}}\leq r < j_{r_{1}}+2$, then $|r_{1}-r|>0$ and if $r=j_{r_{1}+1}$ satisfies (6.1), or if $r=j_{r_{1}}$ then $UWD$ satisfies (4.5).

    \item Now, if $|r_{1}-r|>0$ then \[g(UWD)=U \dots \underbrace{f( w_{l})\dots f(w_{i_{r}+\sigma(a)})}_{DU} \underbrace{w_{i_{r}+\sigma(a)+1}\dots w_{s}}_{UD} \dots D \in S,\]

    \end{enumerate}

    or $r \in [i_{r_{1}+ \sigma(a)},j_{r_{1}+1}]$ with $|r_{1}-1| \geq 0$  satisfies conditions (4.1), (4.2) and (4.3) for $i_{r_{1}+\sigma(a)} \leq r \leq j_{r_{1}+1}$.

\end{enumerate}

Same arguments are used for the cases $r \in [i_{r},j_{r+1- \sigma(a)}]([j_{r}, i_{r+ \sigma(a)}])$ to conclude the lemma. We are done. $\hfill\square$

 \begin{lema}\label{lema3.2}

 \textit{Suppose that $UWD$ is a Dyck path in $S$ and that integers $l$ and $r$ as defined in  $\mathrm{Proposition \hspace{0.1cm}\ref{prop3.2}}$ are such that $l=r$, then the following statements hold}: 

 \begin{enumerate}

     \item [(a)] \textit{If $l \notin \{ j_{s}\}$ then there is an elementary shift to a Dyck path with exactly $n-1$ peaks.}
     \item [(b)] \textit{If $l \in \{ j_{s}\}$ then there is an elementary shift from a Dyck path with exactly $n-1$  peaks to $UWD$.}

 \end{enumerate}

 \end{lema}

\textbf{Proof.} Let $UWD$ be a Dyck path in $S$, and positive integers  $l$ and $r$ with $l=r$. 

\begin{itemize}

    \item [(a)] Suppose $l \notin \{j_{s} \}$ and $l\in [i_{r},j_{r+1-\sigma(a)}]$. If $i_{r}\leq l < j_{r+1-\sigma(a)}$, then  $UWD$ satisfies  (4.1) and (4.2) of Lemma \ref{lema3.1}. In particular, if $l=i_{r}\neq 1$ there is $p'=l-1$ in $[j_{r-\sigma(a)},i_{r}]$ that satisfies the first condition of (5) of Lemma \ref{lema3.1}. The case $l \in [j_{r},i_{r+\sigma(a)}]$ is dual.
    \item [(b)] Suppose $l =j_{r+1-\sigma(a)}$, we have the following cases:

    \begin{enumerate}

        \item [(i)]If $|i_{r}-j_{r+1-\sigma(a)}|=1$ (or $|i_{r+1}-j_{r+1-\sigma(a)}|=1$) and $i_{r}=1$ (or $i_{r+1}=n-1$), then there is $UVD$ which is equal to \[U \underbrace{w_{1}}_{UD}w_{l}\dots D\in S \text{ (or }  U \dots w_{l}\underbrace{w_{n-1}}_{UD}D\in S \text{)},\] and \[Uf(w_{1})w_{l}\dots D=UWD \text{ (or }  U \dots w_{l}f(w_{n-1})D=UWD \text{)}\]

        \item [(ii)]If $|i_{r}-j_{r+1-\sigma(a)}|=1$ (or $|i_{r+1}-j_{r+1-\sigma(a)}|=1$) and $i_{r} \neq 1$ (or $i_{r+1} \neq n-1$) then there is $l'=j_{r-\sigma(a)} $ and $r'=j_{r+1-\sigma(a)}$ (or $l'=j_{r+1-\sigma(a)} $ and $r'=j_{r+2-\sigma(a)}$) such that $UVD$ is equal to \[U \dots \underbrace{w_{l'} \dots w_{r'-1}w_{l}}_{UD} \dots D \text{ (or } U  \dots \underbrace{w_{l}w_{l'+1} \dots w_{r'}}_{UD} \dots D \text{)}\in S\] and \[ U \dots f(w_{l'}) \dots f(w_{r'-1})w_{l} \dots D \text{ (or } U \dots w_{l}f(w_{l'+1})\dots f(w_{r'}) \dots D \text{)}=UWD.\] 

        \item [(iii)]If $|i_{r}-j_{r+1-\sigma(a)}|>1$ (or $|i_{r+1}-j_{r+1-\sigma(a)}|>1$) then there is $UVD$ which  is equal to \[U \dots \underbrace{w_{l-1}}_{UD}w_{l}\dots D \text{ (or } U \dots w_{l}\underbrace{w_{l+1}}_{UD} \dots D\in S \text{)}\] and \[U \dots f(w_{l-1})w_{l} \dots D=UWD \text{ (or } U \dots w_{l}f(w_{l+1}) \dots D=UWD \text{)}.\]  

    \end{enumerate}

\end{itemize}

Similar arguments dually applied can be used to obtain the lemma in the case $l=j_{r}$. We are done. $\hfill\square$

\addtocounter{Nota}{4}

\begin{Nota}\label{Nota3.1}

Note that, in general there is an elementary leftshift and an elementary rightshift  over $S$, and these elementary shifts are disjoint, i.e. if $f_{p_{1}} \circ \dots \circ f_{p_{q}}$ and $f_{p'_{1}} \circ \dots \circ f_{p'_{q'}}$ are elementary left and right shifts, respectively. Then \[ \{p_{1}, \dots ,p_{q} \} \cap \{ p'_{1}, \dots, p'_{q'}\} = \varnothing,\]  these elementary shifts are unique according to Lemma \ref{lema3.1} and Lemma \ref{lema3.2}. If $F^{p}=f_{p_{1}} \circ \dots \circ f_{p_{q}}$ is an elementary leftshift (rightshift) we write $F^{p}_{l}$ ($F^{p}_{r}$). 

\end{Nota}

\addtocounter{prop}{3}

\begin{prop}\label{prop3.3}

\textit{Let $\mathscr{C}=\{ i_{1}, \dots i_{k},j_{1}, \dots j_{m} \}$ be an admissible subchain, then  all Dyck paths of $S$ constitute a connected quiver $Q$}.

\end{prop}

\textbf{Proof.} It suffices to prove that $Q$ is connected, to do that, consider Dyck paths $UWD$ and $UVD$ of $ S$. Then if there is a shift path between $UWD$ and $UVD$ they are connected. Otherwise, Lemmas \ref{lema3.1} and \ref{lema3.2} allow to define a  Dyck path $UW^{(1)}D$  and  a shift path $F^{(1)}=F_{p_{1}}^{(1)} \circ \dots \circ F_{1}^{(1)}$ with $F_{m}^{(1)}=f_{m_{1}}^{(1)} \circ \dots \circ f_{m_{q_{1}}}^{(1)}$ such that \[UWD  \xrightarrow[]{F^{(1)}_{1}}  \dots \xrightarrow[]{F^{(1)}_{p_{1}}} UW^{(1)}D,\] and if there is a shift path from $UVD$ to a $UW^{(1)}D$ then they are connected. If there is not a shift path from $UVD$ to $UW^{(1)}D$,  then there is a Dyck path $UW^{(2)}D$  and  a shift path $F^{(2)}=F_{p_{2}}^{(2)} \circ \dots \circ F_{1}^{(2)}$ with $F_{m}^{(2)}=f_{m_{1}}^{(1)} \circ \dots \circ f_{m_{q_{2}}}^{(2)}$ such that \[UW^{(2)}D  \xrightarrow[]{F^{(2)}_{1}}   \dots  \xrightarrow[]{F^{(2)}_{p_{2}}}  UW^{(1)}D \xleftarrow[]{F^{(1)}_{p_{1}}} \dots \xleftarrow[]{F^{(1)}_{1}} UWD,\] again, if there is a shift path from $UW^{(2)}D$ to a $UVD$ then they are connected.  Since $S$ is finite, the procedure ends in such a way that $UWD$ and $UVD$ are connected and with this argument we are done.  $\hfill\square$ \par\bigskip

Henceforth, we let  $\mathfrak{C_{2n}}$ denote the subcategory of $(\mathfrak{D}_{2n},R_{j_{1} \ldots j_{m}}^{i_{1}\ldots i_{k}})$ whose objects are $k-$linear combinations of Dyck paths of $S$. Lemma \ref{lema3.3} and Proposition \ref{prop3.4} give some properties of the Hom-spaces of this category. 

\addtocounter{lema}{2}

\begin{lema}\label{lema3.3}

\textit{Let $UWD$, $UW'D$, $UW''D$ and $UVD$ be Dyck paths in $\mathfrak{C_{2n}}$ and let $F=F^{1}_{l} \circ F^{2}_{r}$ (resp. $F^{1}_{r} \circ F^{2}_{l}$)  be a shift path $UWD \xrightarrow[]{F^{2}_{r}} UW'D \xrightarrow[]{F^{1}_{l}} UVD$ (resp. $UWD \xrightarrow[]{F^{2}_{l}} UW'D \xrightarrow[]{F^{1}_{r}} UVD$), if there is other shift path $G=G^{1} \circ G^{2}$ such that $UWD  \xrightarrow[]{G^{2}} UW''D \xrightarrow[]{G^{1}} UW''D$ with $UW'D \neq UW''D$ then $G^2=F^{1}_{l}$ and $G^{1}=F^{2}_{r}$ (resp. $G^2=F^{1}_{r}$ and $G^{1}=F^{2}_{l}$).}

\end{lema}

\textbf{Proof.} Let $F=F^{1}_{l} \circ F^{2}_{r}$ be a shift path such that\[U  \dots w_{l_{1}}\dots w_{r_{1}}\dots D \xrightarrow[]{F^{2}_{r}}  U \dots w'_{l_{2}} \dots w'_{r_{2}} \dots D \xrightarrow[]{F^{1}_{l}} U \dots v_{l_{3}}\dots v_{r_{3}}\dots D,\] with $l_{1}=l_{2}$ and $r_{2}=r_{3}$ and suppose that there is other shift path $G=G^{1} \circ G^{2}$ such that \[U  \dots w_{l_{1}}\dots w_{r_{1}}\dots D \xrightarrow[]{G^{2}}  U \dots w''_{l_{4}} \dots w''_{r_{4}} \dots D \xrightarrow[]{G^{1}} U \dots v_{l_{3}}\dots v_{r_{3}}\dots D,\] with $UW'D \neq UW''D$. Given the elementary rightshift $F^{2}_{r}$, then since $G\neq F^{2}_{r}$, it holds that $UWD$  satisfies the conditions of $UW'D$ in order to apply the same elementary leftshift $F^{2}_{l}$, i.e., $F^{2}_{l}=G^{2}$ and $l_{3}=l_{4}$. Since $r_{1}=r_{4}$, $UWD$ and $UW''D$ satisfy the conditions to apply the same elementary rightshift, i.e., $F^{1}_{r}=G^{1}$. Case $F^{1}_{r} \circ F^{2}_{l}$ is obtained via a dual argument. $\hfill\square$

\addtocounter{prop}{1}

\begin{prop}\label{prop3.4}

\textit{If $\mathrm{Hom}_{\mathfrak{C_{2n}}}(UWD,UVD)\neq0$ then $\mathrm{dim}_{k}\hspace{0.1cm}\mathrm{Hom}_{\mathfrak{C_{2n}}}(UWD,UVD)=1$.}

\end{prop}

\textbf{Proof.} Suppose that $\mathrm{Hom}_{\mathfrak{C_{2n}}}(UWD,UVD)\neq0$, then there is a shift path $F$ of the form \[UWD \xrightarrow[]{F^{0}_{x_{0}}} \dots  \xrightarrow[]{F^{i-1}_{x_{i-1}}} UW^{i-1}D \xrightarrow[]{F^{i-1}_{x_{i-1}}}UW^{i}D \xrightarrow[]{F^{i}_{x_{i}}} UW^{i+1}D  \xrightarrow[]{F^{i+1}_{x_{i+1}}}\dots \xrightarrow[]{F^{m}_{x_{m}}}UVD,\]  with $x_{i} \in \{ l, r\}$ and for some $m \in \mathbb{Z}^{+}$. Now, for each pair $F^{i}_{x_{i}}\circ F^{i-1}_{x_{i-1}}$ with $x_{i-1}=l$ and $x_{i}=r$ ($x_{i-1}=r$ and $x_{i}=l$) that satisfies  conditions described in Lemma \ref{lema3.3} there is another shift path $F'$ of the form \[UWD \xrightarrow[]{F^{0}_{x_{0}}} \dots  \xrightarrow[]{F^{i-1}_{x_{i-1}}} UW^{i-1}D \xrightarrow[]{F^{i}_{x_{i}}}UW^{i'}D \xrightarrow[]{F^{i-1}_{x_{i-1}}} UW^{i+1}D  \xrightarrow[]{F^{i+1}_{x_{i+1}}}\dots \xrightarrow[]{F^{m}_{x_{m}}}UVD,\] transforming  $UWD$ and $UVD$. Thus  $F \sim_{R_{j_{1} \ldots j_{m}}^{i_{1}\ldots i_{k}}} F'$. $\hfill\square$

\section{A Categorical Equivalence} \label{section3.2}

In this section, we establish an equivalence between the category $\mathfrak{C_{2n}}$  and the category of representations of a quiver of  Dynkin type $\mathbb{A}_{n}$. 

\subsection{The $\Theta$ Functor}

Given an admissible subchain $\mathscr{C}=\{j_{1},\dots, j_{m},i_{1},\dots, i_{k}\}$, $\mathfrak{C}_{2n}$ the subcategory of $(\mathfrak{D}_{2n},R_{j_{1}\ldots j_{m}}^{i_{1} \ldots i_{k}})$ and $Q$ a quiver of type  $\mathbb{A}_{n-1}$  with $\lbrace i_{1},\dots, i_{k}\rbrace$ and $\lbrace j_{1},\dots, j_{m} \rbrace$ being the sets of sinks and sources, respectively.  Then the $k$-linear additive functor $\Theta: \mathfrak{C}_{2n} \longrightarrow  \mathrm{rep}\hspace{0.1cm} Q$ is defined in such a way that, for an object $UWD\in \mathfrak{C}_{2n}$, it holds that,

\[\Theta (UWD)=(\Theta(w_{i}),\varphi_{\Theta (w_{i},w_{i+1})}),\]

 where  
 
 \begin{equation}
 \Theta (w_{i})=\begin{cases}
  k, & \mbox{if $ w_{i}=UD$,}\\
 0, & \mbox{if  $w_{i}=DU$.}
\end{cases}
\end{equation}  

If  $w_{i},w_{i+1} \in [i_{r},j_{r+1-\sigma(a)}]$ \big ($[j_{r},i_{r+\sigma(a)}]$ \big ) then $s(\Theta (w_{i},w_{i+1}))=i+1$, is the starting point of the corresponding arrow, whereas $t(\Theta (w_{i},w_{i+1}))=i$ is the ending vertex of the corresponding arrow \big ($s(\Theta (w_{i},w_{i+1}))=i$, $t(\Theta (w_{i},w_{i+1}))=i+1$ \big) and,

\[ \varphi_{\Theta (w_{i},w_{i+1})}: \Theta (w_{s(\Theta (w_{i},w_{i+1}))})\longrightarrow \Theta (w_{t(\Theta (w_{i},w_{i+1}))}),\]

\begin{equation}
\varphi_{\Theta (w_{i},w_{i+1})} =\begin{cases}
  1_{k }, & \mbox{if $ w_{i}=UD=w_{i+1}$,}\\
 0, & \mbox{if  $w_{i}=DU$ or $w_{i+1}=DU$.}
\end{cases}
\end{equation}  

Functor $\Theta$ acts on morphisms as follows;\par\bigskip

Let \[f_{q_{2}} \circ \dots \circ  f_{q_{1}}=(1_{1}, \ldots,1_{q_{1}-1}, f_{q_{1}},\ldots, f_{q_{2}},1_{q_{2}+1},\ldots 1_{n-1}), \] be a elementary shift between $UWD$ and $UVD$, then:

 \[\Theta ((1_{1}, \ldots,1_{q_{1}-1}, f_{q_{1}},\ldots, f_{q_{2}},1_{q_{2}+1},\ldots 1_{n-1})),\] \[(\Theta(1_{1}),\ldots,\Theta (1_{q_{1}-1}),\Theta (f_{q_{1}}),\ldots,\Theta (f_{q_{2}}),\Theta (1_{q_{2}+1}),\ldots, \Theta (1_{n-1})),\]

 where  $\Theta(f_{m})=0$  and, 

\begin{equation}
 \Theta (1_{m_{1}})=\begin{cases}
 1_{k}, & \mbox{if $ w_{m_{1}}=UD=v_{m_{1}}$,}\\
0, & \mbox{otherwise,}
\end{cases}
\end{equation} 

for $1\leq m_{1} \leq q_{1}-1$, $ q_{1}\leq m \leq q_{2}$ and $q_{2}+1\leq m_{1} \leq n-1$. 

\addtocounter{Nota}{3}

\begin{Nota}

Note that, it is easy to see that $\Theta$ is an additive covariant functor.

\end{Nota}

\addtocounter{lema}{2}

\begin{lema}\label{lema3.4}

\textit{Let $UWD$ and $UVD$ be Dyck paths of $\mathfrak{C_{2n}}$. If $\mathrm{Hom }_{\mathfrak{C_{2n}}}(UWD,UVD)\neq0$ then $\mathrm{Hom }_{\mathrm{rep}\hspace{0.1cm}Q}(\Theta(UWD),\Theta(UVD))\neq0$}.

\end{lema}

\textbf{Proof.} Suppose $\mathrm{Hom }_{\mathfrak{C_{2n}}}(UWD,UVD)\neq0$, and let $F$ be a shift path $UW^{0}D \xrightarrow[]{F^{0}} UW^{1}D \xrightarrow[]{F^{1}} \dots \xrightarrow[]{F^{m-2}} UW^{m-1}D \xrightarrow[]{F^{m-1}} UW^{m}D$ from $UWD=UW^{0}D$  to $UVD=UW^{m}D$ for some $m \in \mathbb{Z}^{+}$, then there exist $q_{1}$ and $q_{2}$ such that \[ \{q_{1},q_{1}+1, \dots , q_{2}-1 ,q_{2}\}= \bigcap_{i \in J} \text{Supp }UW^{i}D,\]  applying $\Theta$ we obtain the following diagram:

\[\xymatrix{ \dots\ar@{-}[r] &  \Theta(w^{0}_{q_{1}-1}) \ar@{-}[r]^{ \text{\quad} \text{\quad} a^{0}_{q_{1}-1}} \ar[d]_{c^{0}_{q_{1}-1}}&   \ar@{-}[r]^{1} k \ar[d]^{1} & \dots\ar@{-}[r]^{1} & k  \ar[d]_{1}& \Theta(w^{0}_{q_{2}+1})\ar@{-}[r] \ar@{-}[l]_{a^{0}_{q_{2}} \text{\quad} \text{\quad} } \ar[d]^{d^{0}_{q_{2}+1}}& \dots \\ 
\dots \ar@{-}[r] & \Theta(w^{1}_{q_{1}-1}) \ar@{-}[r]^{ \text{\quad} \text{\quad} a^{1}_{q_{1}-1}} \ar[d]_{c^{1}_{q_{1}-1}}& k \ar[d]^{1} \ar@{-}[r]^{1}& \dots \ar@{-}[r]^{1} & k  \ar[d]_{1}& \Theta(w^{1}_{q_{2}+1})\ar@{-}[r] \ar@{-}[l]_{a^{1}_{q_{2}}\text{\quad} \text{\quad} } \ar[d]^{d^{1}_{q_{2}+1}}& \dots \\& \vdots  \ar[d]_{c^{m-2}_{q_{1}-1}}& \vdots \ar[d]^{1}&&\vdots \ar[d]_{1}& \vdots \ar[d]^{d^{m-2}_{q_{2}+1}} &\\ \dots \ar@{-}[r]& \Theta(w^{m-1}_{q_{1}-1})\ar@{-}[r]^{ \text{\quad} \text{\quad} a^{m-1}_{q_{1}-1}} \ar[d]_{c^{m-1}_{q_{1}-1}}& k \ar@{-}[r]^{1} \ar[d]^{1} & \dots \ar@{-}[r]^{1}& k  \ar[d]_{1}& \Theta(w^{m-1}_{q_{2}+1}) \ar@{-}[r] \ar@{-}[l]_{a^{m-1}_{q_{2}}\text{\quad} \text{\quad} } \ar[d]^{d^{m-1}_{q_{2}+1}}& \dots \\ \dots \ar@{-}[r] & \Theta(w^{m}_{q_{1}-1}) \ar@{-}[r]^{ \text{\quad} \text{\quad} a^{m+1}_{q_{1}-1}}& k \ar@{-}[r]^{1}& \dots\ar@{-}[r]^{1} & k & \Theta(w^{m}_{q_{2}+1})\ar@{-}[r] \ar@{-}[l]_{a^{m}_{q_{2}}\text{\quad} \text{\quad} }& \dots }\]
\[\text{Diagram 1.}\]

where $c^{i}_{q_{1}}-1$, $a^{i}_{q_{1}-1}, a^{i}_{q_{2}}, d^{i}_{q_{2}+1}\in\{0,k\}$, squares in the diagram are commutative between $q_{1}$ and $q_{2}$ (independently of the chosen orientation). For the sub-shift path $F^{(x, y)}$ to $F$ with $0 \leq x \leq y \leq m-1$   there exist positive integers $q_{1}^{(x, y)}$ and $q_{2}^{(x, y)}$ such that \[S^{(x, y)}= \{q_{1}^{(x, y)},q_{1}^{(x, y)}+1, \dots, q_{2}^{(x, y)}-1, q_{2}^{(x, y)} \} =\bigcap_{i \in J^{(x ,y)}} \text{Supp }UW^{i}D,\]  and for the diagrams

\[ \xymatrix{ \Theta(w^{x}_{q^{(x, y)}_{1}-1}) \ar@{-}[r]^{ \text{\quad} \text{\quad} a^{x}_{q^{(x, y)}_{1}-1}} \ar[d]_{0}&    k \ar[d]^{1} \\ \Theta(w^{y}_{q^{(x, y)}_{1}-1}) \ar@{-}[r]^{ \text{\quad} \text{\quad} a^{y}_{q^{(x, y)}_{1}-1}}&k} \text { and } \xymatrix{k  \ar[d]_{1}& \Theta(w^{x}_{q^{(x ,y)}_{2}+1}) \ar@{-}[l]_{a^{x}_{q^{(x, y)}_{2}} \text{\quad} \text{\quad} } \ar[d]^{0}\\ k & \Theta(w^{y}_{q^{(x, y)}_{2}+1}) \ar@{-}[l]_{a^{y}_{q^{(x, y)}_{2}}\text{\quad} \text{\quad} }}\]  \[ \text{Diagram 2.  \qquad } \text{ \quad } \text{ \qquad  Diagram 3.} \]

we have the following cases:

 \begin{enumerate}[(1)]

     \item If $q^{(x, y)}_{1} \in [i_{r},j_{r+1-\sigma(a)}]$ ($i_{r}< q^{(x ,y)}_{1} \leq j_{r+1-\sigma(a)} $) four cases must be considered.

     \begin{enumerate}

          \item [(1.1)] If $\Theta(w^{x}_{q^{(x, y)}_{1}-1})=\textit{k}$ and $\Theta(w^{y}_{q^{(x, y)}_{1}-1})=\textit{k}$, $q_{1}^{(x, y)}$ belong to $ S^{(x, y)}$, which is a contradiction.

          \item [(1.2)] If $\Theta(w^{x}_{q^{(x, y)}_{1}-1})=\textit{k}$ and $\Theta(w^{y}_{q^{(x, y)}_{1}-1})=0$, then the Diagram  2 commutes.

          \item  [(1.3)] If $\Theta(w^{x}_{q^{(x, y)}_{1}-1})=0$ and $\Theta(w^{y}_{q^{(x, y)}_{1}-1})=\textit{k}$, then  there is an elementary shift $f_{q^{(x, y)}_{1}-1}$ on the interval and this is again a contradiction.

          \item [(1.4)]  If $\Theta(w^{x}_{q^{(x, y)}_{1}-1})=0$ and $\Theta(w^{y}_{q^{(x, y)}_{1}-1})=0$, then the Diagram 2 commutes.

     \end{enumerate}

       \item If $q^{(x, y)}_{1} \in [j_{r+1-\sigma(a)},i_{r+1}]$ ($j_{r+1-\sigma(a)} <  q^{(x, y)}_{1} \leq  i_{r+1}$), the conditions (1.1)-(1.4) are satisfied on the interval.

       \begin{enumerate}

           \item [(2.1)] If $\Theta(w^{x}_{q^{(x ,y)}_{1}-1})=\textit{k}$ and $\Theta(w^{y}_{q^{(x, y)}_{1}-1})=0$, then they satisfy condition (1.3).

           \item [(2.2)]  If $\Theta(w^{x}_{q^{(x, y)}_{1}-1})=0$ and $\Theta(w^{y}_{q^{(x, y)}_{1}-1})=\textit{k}$, then they satisfy condition (1.2).

       \end{enumerate}

       \item Case $q^{(x, y)}_{2} \in [i_{r},j_{r+1-\sigma(a)}]$ is similar to case (2) for the Diagram 3.

       \item  Case $q^{(x ,y)}_{2} \in [j_{r+1-\sigma(a)},i_{r+1}]$ is similar to case (1) for the Diagram 3.

 \end{enumerate}

therefore the Diagram 1 commutes. Since the cases over $[j_{r},i_{r+\sigma(a)}]$ can be showed by using dual arguments. We are done. $\hfill\square$

\begin{lema}\label{lema3.5}

\textit{Functor $\Theta$ is faithful and full.}

\end{lema}

\textbf{Proof.} Let $\phi$ be the map \[\phi: \mathrm{Hom }_{\mathfrak{C_{2n}}}(UWD,UVD) \rightarrow \text{Hom }_{\mathrm{rep}\hspace{0.1cm}Q}(\Theta(UWD),\Theta(UVD)),\] such that $\phi(\lambda F)= \lambda \Theta (F)$ with $F=(1_{1}, \dots ,1_{q_{1}-1}, f_{q_{1}}, \dots f_{q_{2}},1_{q_{2}+1}, \dots, 1_{n-1})$, for some $1 \leq q_{1},q_{2} \leq n-1$ and $\lambda \in k$. Note, $\phi$ is well defined and Lemma \ref{lema3.4} allows us to observe that the image of a non-zero morphism in $\mathfrak{C_{2n}}$ is a non-zero morphism in $\mathrm{rep}\hspace{0.1cm} Q$. Thus, $\phi$ is surjective and injective. $\hfill\square$

\addtocounter{teor}{12}

\begin{teor}\label{teor3.1}

\textit{Functor $\Theta$ is a categorical equivalence between the categories $\mathfrak{C}_{2n}$ and $\mathrm{rep}\hspace{0.1cm} Q$.}

\end{teor}

\textbf{Proof.} Lemma \ref{lema3.5} implies that functor $\Theta$ is faithful and full. Now, let $(M_{i},\varphi_{\alpha})_{i \in Q_{0}, \alpha \in Q_{1}}$ be an indecomposable representation in $\mathrm{rep}\hspace{0.1cm}Q$ of the form

\begin{center}

\begin{picture}(315,15)

\put(0,0){$\xymatrix{0 \ar@{-}[r]& \cdots \ar@{-}[r]& k \ar@{-}[r]^{\text{  }1}& k  \ar@{-}[r]^{1}& \cdots  \ar@{-}[r]^{1}& k  \ar@{-}[r]^{1 \text{  }}& k \ar@{-}[r] & \cdots \ar@{-}[r]&0}$}
\put(75,10){$\overbrace{}^{q_1}$}

\put(220,10){$\overbrace{}^{q_2}$}

\end{picture}

\end{center}

with  $\lbrace i_{1},\dots, i_{k}\rbrace$ and $\lbrace j_{1},\dots, j_{m} \rbrace$  the sets of sinks and sources respectively. Let $\varphi_{1}: \{ 0,k \} \rightarrow \{ DU,UD\}$ be a map  such that $\varphi_{1}(k)=UD$ and $\varphi_{1}(0)=DU$. Define the Dyck path $UWD$ such that \[UWD=U\overbrace{w_{1}\dots w_{q_{1}-1}}^{DU} \overbrace{w_{q_{1}}\dots w_{q_{2}}}^{UD} \overbrace{w_{q_{2}+1}\dots w_{n-1}}^{DU}D. \]  Proposition \ref{prop3.2} allows us to observe that  $UWD$ has $n-1$ peaks over  $\{j_{1},\dots, j_{m},i_{1},\dots, i_{k}\}$ and $\Theta(UWD)=(M_{i},\varphi_{\alpha})_{i \in Q_{0}, \alpha \in Q_{1}}$. Thus,  $\Theta$ is essentially surjective. $\hfill\square$

\addtocounter{corol}{15}

\begin{corol}\label{corol4.1}

\textit{There exists a bijection $\varphi$ between the set of representatives of indecomposable representations of $\mathrm{rep}\hspace{0.1cm}Q$ (denoted  $\mathrm {Ind}(\mathrm{rep}\hspace{0.1cm}Q)$) and the set of Dyck paths of length $2n$ with exactly $n-1$ peaks. }

\end{corol}

\textbf{Proof.} The Narayana number with  exactly $n-1$ peaks over all Dyck paths of length $2n$ is the triangular number $T_{n-1}=\frac{(n-1)(n)}{2}$, which is equal to the number of indecomposable representations of $\mathrm{rep}\hspace{0.1cm}Q$, then we define $\varphi: S \rightarrow \text {Ind }(\mathrm{rep}\hspace{0.1cm}Q)$ such that $\varphi(UWD)=\Theta (UWD)$. $\hfill\square$

\begin{corol}\label{corol4.2}

\textit{The category $\mathfrak{C}_{2n}$ is an abelian category.}

\end{corol}

\subsection{Properties of the Category $\mathfrak{C}_{2n}$}

In this section, we introduce some properties of $\mathfrak{C}_{2n}$ regarding simple, projective and injective indecomposable  objects, also  we construct the Auslander-Reiten quiver for algebras of Dynkin type $\mathbb{A}_{n-1}$. Some conditions for morphisms between objects of the category are introduced as well.

\addtocounter{teor}{2}

\begin{teor}\label{teor3.2}

\textit{Let $\mathscr{C}=\{j_{1},\dots, j_{m},i_{1},\dots, i_{k}\}$ be an admissible subchain, and  let $\mathfrak{C}_{2n}$ be the corresponding category, then} 

\begin{enumerate}[(i)]

\item \textit{Indecomposable simple objects of $\mathfrak{C}_{2n}$ are objects of the form} \[S(x)=US(w^{x}_{1})\ldots S(w^{x}_{n})D\] \textit{where} 

\begin{equation}
 S(w_{y}^{x})=\begin{cases}
  UD, & \mbox{if $ x=y$,}\\
DU, & {otherwise.}
\end{cases}
\end{equation}

\item \textit{Indecomposable projective objects of $\mathfrak{C}_{2n}$ have the form $P(x)=UP(w^{x}_{1})\ldots P(w^{x}_{n})D$ where} 

\begin{equation}
P(w_{x}^{y})=\begin{cases}
  UD, & \mbox{if $ x, y \in [i_{r},j_{r+1-\sigma(a)}]$ $([j_{r},i_{r+\sigma(a)}) $ and $y\leq x$  $(x\leq y),$ }\\
DU, & {otherwise.}
\end{cases}
\end{equation}

\item \textit{Indecomposable injective objects of $\mathfrak{C}_{2n}$ have the form $I(i)=UI(w^{x}_{1})\ldots I(w^{x}_{n})D$ where}

\begin{equation}
I(w_{x}^{y})=\begin{cases}
  UD, & \mbox{if $ x, y \in [i_{r},j_{r +1-\sigma(a)}]$ $([j_{r},i_{r+\sigma(a)}])$ and $x\leq y$ $(y \leq x),$ }\\
DU, & {otherwise.}
\end{cases}
\end{equation}

\end{enumerate}

\end{teor}

\textbf{Proof.} (i) Let $S(x)=(S(x)_{y},\varphi_{\alpha})$ be an indecomposable simple object of $\mathrm{rep}\hspace{0.1cm}Q$ such that $S(x)_{y}=k$ if $x=y$ and $S(x)_{y}=0$ if $x\neq y$. Functor $\Theta$ allows us to observe that, there is $UWD\in\mathfrak{C}_{2n}$ satisfying the required conditions. \par\bigskip

(ii) Let $P(x)=(P(x)_{y},\varphi_{\alpha})$ be an indecomposable projective object of $\mathrm{rep}\hspace{0.1cm}Q$, if $P(x)_{y}=k$ then there is a path from $x$ to $y$, as well as, a source $j_{r+1-\sigma(a)}$ ($j_{r}$) and a sink $i_{r}$ ($i_{r+\sigma(a)}$) such that $i_{r}\leq y \leq x \leq j_{r+1-\sigma(a)}$ ($j_{r} \leq x \leq y \leq i_{r+\sigma(a)}$), and $P(x)_{y}=0$. Thus, there is not a path between $x$ and $y$, then functor $\Theta$  determines an object $UWD$ of $\mathfrak{C}_{2n}$ with $i_{1},\dots i_{k},j_{1},\dots j_{m}$ being an admissible subchain satisfying the required conditions. Case (iii) follows by dually applying the arguments used in the case (ii). $\hfill\square$

\addtocounter{corol}{1}

\begin{corol}\label{corol3.3}

\textit{The indecomposable simple objects of $\mathfrak{C}_{2n}$ have exactly a subsequence $UUDD$.}

\end{corol}

 \textbf{Proof.} Let $S(x)$ be an indecomposable simple object of $\mathfrak{C}_{2n}$, then the identity \[S(x)=U \dots S(w_{x-1}^{x})S(w_{x}^{x})S(w_{x+1}^{x})\dots D=U \dots DU \dots \underbrace{DU}_{x-1} \underbrace{UD}_{x} \underbrace{DU}_{x+1} \dots DU \dots D\] has place as a consequence of Theorem  \ref{teor3.2}. $\hfill\square$ 

\addtocounter{Nota}{7}

\begin{Nota}\label{Nota3.4}

 The Auslander-Reiten translate can be obtained by using the Coxeter transformation and the dimension vector associated to the support of a Dyck path in $\mathfrak{C}_{2n}$.

\end{Nota}

Morphisms in $\mathfrak{C}_{2n}$ also have the following properties.\par\bigskip

Let $UWD$  be a Dyck path of $\mathfrak{C}_{2n}$, then

\begin{itemize}

    \item $p_{UWD}=w_{t}$ and $b_{UWD}=\text{max } \{ w_{s} \text{ }| \text{ } w_{i_{r}} \leq w_{s} \leq w_{j_{r+1-\sigma(a)}} \text{, } w_{s}=UD \}$ over $[i_{r},j_{r+1-\sigma(a)}]$,

    \item $p^{UWD}=\text{min } \{ w_{s} \text{ }| \text{ } w_{j_{r}} \leq w_{s} \leq w_{i_{r+\sigma(a)}} \text{, } w_{s}=UD \}$ and $b^{UWD}=w_{t}$ over $[j_{r},i_{r+\sigma(a)}].$

\end{itemize}

\par\bigskip

\addtocounter{teor}{2}

\begin{teor}\label{teor3.3}

\textit{The vector space $\mathrm{Hom}_{\mathfrak{C}_{2n}}(UWD,UVD)\neq0$ if and only}  

\begin{enumerate}[(i)]

\item  $\text{Supp}(UWD)\cap \text{Supp} (UVD)\neq \varnothing$,
\item   $p_{UWD}\leq p_{UVD}$ and $b_{UWD}\leq b_{UVD}$ \textit{over} $[i_{r},j_{r+1-\sigma(a)}]$,
\item  $p^{UWD}\geq p^{UVD}$ and $b^{UWD}\geq b^{UVD}$ \textit{over} $[j_{r},i_{r+\sigma(a)}]$,

\end{enumerate}

\textit{for all $[i_{r},j_{r+1-\sigma(a)}]$, $[j_{r},i_{r+\sigma(a)}]$ such that $i_{r} \leq q \leq j_{r+1-\sigma(a)}$ and $j_{r} \leq q \leq i_{r+\sigma(a)}$  with $q \in \text{Supp}(UWD) \cap \text{Supp} (UVD) $. }

\end{teor}

\textbf{Proof.} The result follows as a consequence of the definition of the functor $\Theta$ and the construction of Lemma 3.1. $\hfill\square$ \\

\begin{figure}[h!]

\begin{center}

\begin{picture}(160,15)

\put(3,0){$Q=$}

\put(30,3){\circle{5}} \put(60,3){\circle{5}} \put(90,3){\circle{5}} \put(120,3){\circle{5}} \put(150,3){\circle{5}}
\put(62,3) {\vector(1,0){25}}
\put(57.7,3){\vector(-1,0){25}}
\put(122,3) {\vector(1,0){25}}
\put(117.7,3){\vector(-1,0){25}}

\put(28,8){$_{1}$}
\put(58,8){$_{2}$}
\put(88,8){$_{3}$}
\put(118,8){$_{4}$}
\put(148,8){$_{5}$}

\end{picture}

\end{center}

\bigskip
\bigskip

\begin{center}

\begin{picture}(276,228)

\multiput(0,0)(96,0){3} {\line(1,1){36}}
\multiput(48,48)(96,0){3} {\line(1,1){36}}
\multiput(0,96)(96,0){3} {\line(1,1){36}}
\multiput(48,144)(96,0){3} {\line(1,1){36}}
\multiput(0,192)(96,0){3} {\line(1,1){36}}
 
\multiput(0,0)(6,6){4} {\line(0,1){6}}
\multiput(0,6)(6,6){4} {\line(1,0){6}} 
\put(24,24){\line(0,1){12}} 
\put(24,36){\line(1,0){12}} 

\multiput(96,0)(6,6){2} {\line(0,1){6}}
\multiput(96,6)(6,6){2} {\line(1,0){6}}
\put(108,12){\line(0,1){12}} 
\put(108,24){\line(1,0){6}} 
\put(114,24){\line(0,1){6}}
\put(114,30){\line(1,0){12}} 
\put(126,30){\line(0,1){6}}
\put(126,36){\line(1,0){6}}

\put(192,0){\line(0,1){12}}
\put(192,12){\line(1,0){6}}
\put(198,12){\line(0,1){6}}
\put(198,18){\line(1,0){12}}
\multiput(210,18)(6,6){3} {\line(0,1){6}}  
\multiput(210,24)(6,6){3} {\line(1,0){6}}

\multiput(48,48)(6,6){2} {\line(0,1){6}} 
\multiput(48,54)(6,6){2} {\line(1,0){6}}
\put(60,60){\line(0,1){12}} 
\multiput(60,72)(6,6){2} {\line(1,0){6}}
\multiput(66,72)(6,6){2} {\line(0,1){6}}
\put(72,84){\line(1,0){12}} 

\put(144,48){\line(0,1){12}} 
\multiput(144,60)(6,6){3} {\line(1,0){6}}
\multiput(150,60)(6,6){3} {\line(0,1){6}} 
\put(162,78){\line(1,0){12}}
\put(174,78){\line(0,1){6}}
\put(174,84){\line(1,0){6}}

\put(240,48){\line(0,1){6}}
\put(240,54){\line(1,0){6}}
\put(246,54){\line(0,1){12}}
\put(246,66){\line(1,0){12}} 
\multiput(258,66)(6,6){3} {\line(0,1){6}}
\multiput(258,72)(6,6){3} {\line(1,0){6}}

\multiput(0,96)(6,6){2} {\line(0,1){6}}
\multiput(0,102)(6,6){2} {\line(1,0){6}}
\put(12,108){\line(0,1){12}}
\put(12,120){\line(1,0){12}}
\multiput(24,120)(6,6){2} {\line(0,1){6}}
\multiput(24,126)(6,6){2} {\line(1,0){6}}

\put(96,96){\line(0,1){12}}
\multiput(96,108)(6,6){4} {\line(1,0){6}}
\multiput(102,108)(6,6){4} {\line(0,1){6}}
\put(120,132){\line(1,0){12}}

\multiput(192,96)(30,30){2} {\line(0,1){6}}
\multiput(192,102)(30,30){2} {\line(1,0){6}}
\put(198,102){\line(0,1){12}}
\multiput(198,114)(6,6){2} {\line(1,0){6}}
\multiput(204,114)(6,6){2} {\line(0,1){6}}
\put(210,126){\line(1,0){12}}

\put(48,144){\line(0,1){12}}
\multiput(48,156)(6,6){2} {\line(1,0){6}}
\multiput(54,156)(6,6){2} {\line(0,1){6}}
\put(60,168){\line(1,0){12}}
\multiput(72,168)(6,6){2} {\line(0,1){6}}
\multiput(72,174)(6,6){2} {\line(1,0){6}}

\put(144,144){\line(0,1){6}}
\put(144,150){\line(1,0){6}}
\put(150,150){\line(0,1){12}}
\multiput(150,162)(6,6){3} {\line(1,0){6}}
\multiput(156,162)(6,6){3} {\line(0,1){6}}
\put(168,180){\line(1,0){12}}

\multiput(240,144)(6,6){3} {\line(0,1){6}}
\multiput(240,150)(6,6){3} {\line(1,0){6}}
\put(258,162){\line(0,1){12}}
\put(258,174){\line(1,0){12}}
\put(270,174){\line(0,1){6}}
\put(270,180){\line(1,0){6}}

\put(0,192){\line(0,1){12}}
\put(0,204){\line(1,0){12}}
\multiput(12,204)(6,6){4} {\line(0,1){6}}
\multiput(12,210)(6,6){4} {\line(1,0){6}}

\put(96,192){\line(0,1){6}}
\put(96,198){\line(1,0){6}}
\put(102,198){\line(0,1){12}}
\put(102,210){\line(1,0){6}}
\put(108,210){\line(0,1){6}}
\put(108,216){\line(1,0){12}}
\multiput(120,216)(6,6){2} {\line(0,1){6}}
\multiput(120,222)(6,6){2} {\line(1,0){6}}

\multiput(192,192)(6,6){3} {\line(0,1){6}}
\multiput(192,198)(6,6){3} {\line(1,0){6}}
\put(210,210){\line(0,1){12}}
\put(210,222){\line(1,0){6}}
\put(216,222){\line(0,1){6}}
\put(216,228){\line(1,0){12}}

\multiput(37,37)(96,0){3} {\vector(1,1){10}}
\multiput(85,85)(96,0){2} {\vector(1,1){10}}
\multiput(37,133)(96,0){3} {\vector(1,1){10}}
\multiput(85,181)(96,0){2} {\vector(1,1){10}}
\multiput(72,60)(96,0){2} {\vector(1,-1){30}}
\multiput(24,108)(96,0){3} {\vector(1,-1){30}}
\multiput(72,156)(96,0){2} {\vector(1,-1){30}}
\multiput(24,204)(96,0){3} {\vector(1,-1){30}}

\end{picture}

    \end{center}

\caption{ Quiver $Q$ and the Auslander-Reiten quiver of $\mathrm{rep}\hspace{0.1cm}Q$.}\label{figure3.3}

\end{figure}
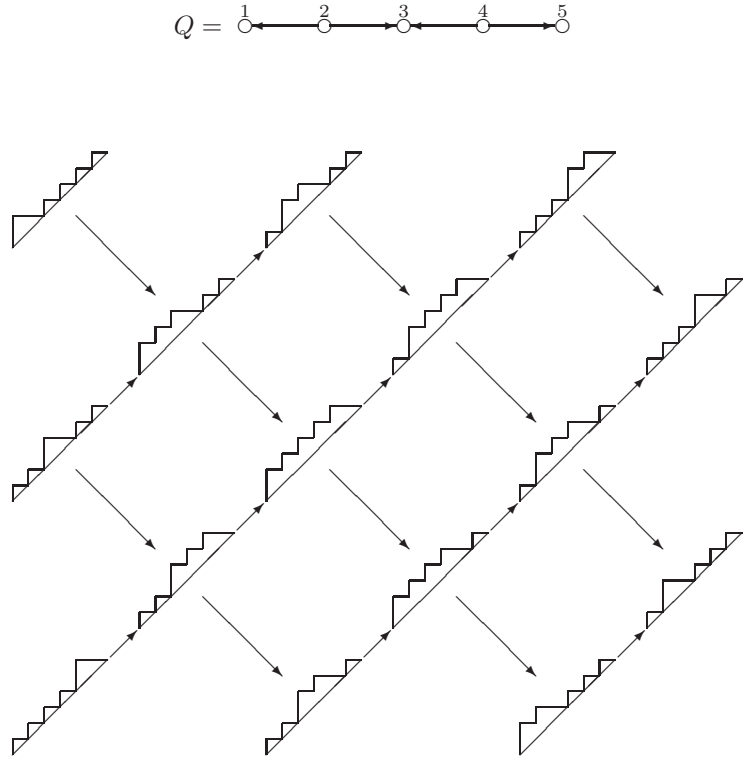

Figure \ref{figure3.3} describes a quiver $Q$ of type $\mathbb{A}_{5}$ and the Auslander-Reiten quiver of $\mathrm{rep}\hspace{0.1cm}Q$.\\

\subsection{A Relationship with Some Nakayama Algebras}

In \cite{Marci} Marczinzik, Rubey and Stump presented a connection between the Auslander-Reiten quiver of Nakayama algebras and Dyck paths. In such a work for a Nakayama algebra $\mathcal{A}$,  they associated the vector space dimension of the indecomposable projective modules $e_{i}\mathcal{A}$ to a Dyck path, this vector is called the Kupisch series. If we take a Nakayama algebra $\mathcal{A}=kQ/I$, with $I= \langle x_{3}x_{4},x_{1}x_{2}x_{3} \rangle$, then the Kupisch series of $kQ/I$ is $[3,3,2,2,1]$, and the Auslander-Reiten quiver of $kQ/I$ has the shape described in Figure \ref{figure3.5}.

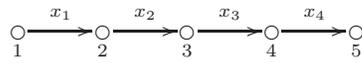
\begin{figure}[h!]

 \begin{center}

 \begin{picture}(140,25)

\put(0,8){{{{\xymatrix @R=0.5cm{ \begin{picture}(2,0)\put(1,2){\circle{5}} \end{picture} \ar@{->}[r] &\begin{picture}(2,0)\put(1,2){\circle{5}} \end{picture}  \ar@{->}[r]&  \begin{picture}(2,0)\put(1,2){\circle{5}} \end{picture} \ar@{->}[r]& \begin{picture}(2,0)\put(1,2){\circle{5}} \end{picture} \ar@{->}[r] & \begin{picture}(2,0)\put(1,2){\circle{5}} \end{picture}}}}}}

\put(16,15){$_{x_{1}}$}
\put(48,15){$_{x_{2}}$}
\put(80,15){$_{x_{3}}$}
\put(112,15){$_{x_{4}}$}
\put(2,0){$_{1}$}
\put(34,0){$_{2}$}
\put(66,0){$_{3}$}
\put(98,0){$_{4}$}
\put(130,0){$_{5}$}

\end{picture}

\end{center}

\caption{Quiver $Q$ of type $\mathbb{A}_{5}$.}\label{figure3.4}

\end{figure}

 \vspace{0.5cm} 

\begin{figure}[h!]

\begin{center}
\begin{picture}(56,56)

\put(0,0){\line(1,1){56}}
\put(0,0){\line(0,1){14}}
\put(0,14){\line(1,0){14}}
\put(14,14){\line(0,1){28}}
\put(14,42){\line(1,0){14}}
\put(28,42){\line(0,1){14}}
\put(28,56){\line(1,0){28}}
\multiput(14,28)(4.3,0){3} {\line(1,0){3.3}} 
\multiput(28,42)(4.3,0){3} {\line(1,0){3.3}} 
\multiput(28,28)(0,4.3){3} {\line(0,1){3.3}} 
\multiput(42,42)(0,4.3){3} {\line(0,1){3.3}} 

\end{picture}
\end{center}

\caption{ Dyck path associated to $kQ/I$.}\label{figure3.5}

\end{figure}
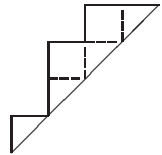 

Let $\mathfrak{C}_{2(n+1)}$  be the category with the admissible subchain $1<n$, $j_{1}=1$ and $i_{1}=n$, and let $D_{i}$ be the sets

\begin{equation}
\begin{array}{l } 
      D_{1}= \{ X \in Ob(\mathfrak{C}_{2(n+1)}) \text{ }| \text{ }  w_{1}=UD \}, \\
      \\
   D_{i}= \{ X \in Ob(\mathfrak{C}_{2(n+1)}) \text{ }| \text{ }  w_{m}=DU, \text{ }1\leq m \leq i-1 \},
   \end{array}
\end{equation}

for $1 < i \leq n$. Then, we take the subset $D_{i,j}  \subseteq D_{i}$, 

\begin{equation}
D_{i,j_{i}}=\{ Y \in D_{i}  \text{ }| \text{ } i \leq r_{Y} \leq m(i,j_{i})+i-1 \},
\end{equation}

such that the vector $v=(n-(m(i,j_{i})+i-1))_{i=1}^{n}$ constitutes an integer partition with $n$ parts.  Now, let  $\mathfrak{N}_{v}$ be the  full subcategory of  $\mathfrak{C}_{2(n+1)}$ whose objects  are $k-$linear combinations of the Dyck paths in the following set 

\begin{equation}\label{equation3.3}
\mathscr{L}= \bigcup_{i=1}^{n}D_{i,j_{i}},
\end{equation}

 and morphisms defined by the category $\mathfrak{C}_{2n(n+1)}$. \par\bigskip

 We assume the following numbering and orientation for a quiver $Q$ associated to a Nakayama algebra 

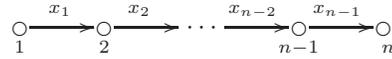
\begin{figure}[h!]

 \begin{center}

 \begin{picture}(148,30)

\put(0,8){{{{\xymatrix @R=0.5cm{ \begin{picture}(2,0)\put(1,2){\circle{5}} \end{picture} \ar@{->}[r] &\begin{picture}(2,0)\put(1,2){\circle{5}} \end{picture}  \ar@{->}[r]& \cdots \ar@{->}[r]& \begin{picture}(2,0)\put(1,2){\circle{5}} \end{picture} \ar@{->}[r] & \begin{picture}(2,0)\put(1,2){\circle{5}} \end{picture} }}}}}

\put(2,0){$_{1}$}
\put(34,0){$_{2}$}
\put(102,0){$_{n-1}$}
\put(141,0){$_{n}$}
\put(15,15){$_{x_{1}}$}
\put(45,15){$_{x_{2}}$}
\put(83,15){$_{x_{n-2}}$}
\put(115,15){$_{x_{n-1}}$}

\end{picture}

\end{center}

\caption{Quiver $Q$ of type $\mathbb{A}_{n}$.}\label{figure3.6}

\end{figure}

The functor $\Theta'$ between the category $\mathfrak{N}_{v}$ and the category of representations of $kQ/I$ where $kQ/I$ is a Nakayama algebra with Kupisch series $[m(1, j_{1}),\dots, m(n, j_{n})]$ is defined in such a way that, $\Theta'(UWD)=\Theta(UWD)$ and $\Theta'(F)=\Theta(F)$ for $UWD\in\mathscr{L}$ and $F$ being an elementary shift in $\mathfrak{N}_{v}$.

\addtocounter{corol}{2}

\begin{corol}

The functor $\Theta'$ is an equivalence of categories.

\end{corol}

\textbf{Proof.} It is a direct consequence of Theorem \ref{teor3.1}.$\hfill\square$\\

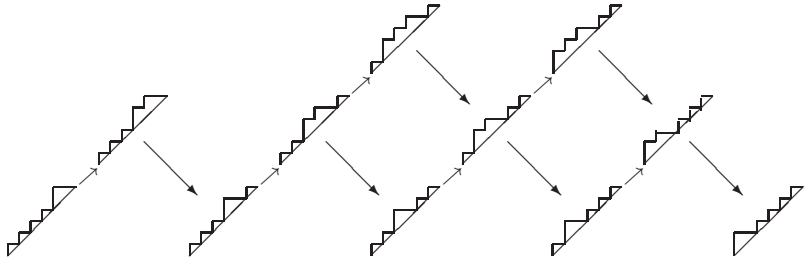
\begin{figure}[h!]
\begin{center}
\begin{picture}(301,70)
\multiput(0,0)(68.8,0){5} {\line(1,1){25.8}}
\multiput(34.4,34.4)(68.8,0){4} {\line(1,1){25.8}}
\multiput(137.6,68.8)(68.8,0){2} {\line(1,1){25.8}}

\multiput(0,0)(4.3,4.3){4} {\line(0,1){4.3}}
\multiput(0,4.3)(4.3,4.3){4} {\line(1,0){4.3}}
\put(17.2,17.2){\line(0,1){8.6}}
\put(17.2,25.8){\line(1,0){8.6}}

\multiput(68.8,0)(4.3,4.3){3} {\line(0,1){4.3}}
\multiput(68.8,4.3)(4.3,4.3){3} {\line(1,0){4.3}}
\put(81.7,12.9){\line(0,1){8.6}}
\put(81.7,21.5){\line(1,0){8.6}}
\put(90.3,21.5){\line(0,1){4.3}}
\put(90.3,25.8){\line(1,0){4.3}}

\multiput(137.6,0)(4.3,4.3){2} {\line(0,1){4.3}}
\multiput(137.6,4.3)(4.3,4.3){2} {\line(1,0){4.3}}
\put(146.2,8.6){\line(0,1){8.6}}
\put(146.2,17.2){\line(1,0){8.6}}
\multiput(154.8,17.2)(4.3,4.3){2} {\line(0,1){4.3}}
\multiput(154.8,21.5)(4.3,4.3){2} {\line(1,0){4.3}}

\put(206,0){\line(0,1){4.3}}
\put(206,4.3){\line(1,0){4.3}}
\put(210.7,4.3){\line(0,1){8.6}}
\put(210.7,12.9){\line(1,0){8.6}}
\multiput(219.3,12.9)(4.3,4.3){3} {\line(0,1){4.3}}
\multiput(219.3,17.2)(4.3,4.3){3} {\line(1,0){4.3}}

\put(274.8,0){\line(0,1){8.6}}
\put(274.8,8.6){\line(1,0){8.6}}
\multiput(283.4,8.6)(4.3,4.3){4} {\line(0,1){4.3}}
\multiput(283.4,12.9)(4.3,4.3){4} {\line(1,0){4.3}}

\multiput(34.4,34.4)(4.3,4.3){3} {\line(0,1){4.3}}
\multiput(34.4,38.7)(4.3,4.3){3} {\line(1,0){4.3}}
\put(47.3,47.3){\line(0,1){8.6}}
\put(47.3,55.9){\line(1,0){4.3}}
\put(51.6,55.9){\line(0,1){4.3}}
\put(51.6,60.2){\line(1,0){8.6}}

\multiput(103.2,34.4)(4.3,4.3){2} {\line(0,1){4.3}}
\multiput(103.2,38.7)(4.3,4.3){2} {\line(1,0){4.3}}
\put(111.8,43){\line(0,1){8.6}}
\put(111.8,51.6){\line(1,0){4.3}}
\put(116.1,51.6){\line(0,1){4.3}}
\put(116.1,55.9){\line(1,0){8.6}}
\put(124.7,55.9){\line(0,1){4.3}}
\put(124.7,60.2){\line(1,0){4.3}}

\put(172,34.4){\line(0,1){4.3}}
\put(172,38.7){\line(1,0){4.3}}
\put(176.3,38.7){\line(0,1){8.6}}
\put(176.3,47.3){\line(1,0){4.3}}
\put(180.6,47.3){\line(0,1){4.3}}
\put(180.6,51.6){\line(1,0){8.6}}
\multiput(189.2,51.6)(4.3,4.3){2} {\line(0,1){4.3}}
\multiput(189.2,55.9)(4.3,4.3){2} {\line(1,0){4.3}}

\put(240.8,34.4){\line(0,1){8.6}}
\put(240.8,43){\line(1,0){4.3}}
\put(245.1,43){\line(0,1){4.3}}
\put(245.1,46.3){\line(1,0){8.6}}
\multiput(253.7,46.3)(4.3,4.3){3} {\line(0,1){4.3}}
\multiput(253.7,51.6)(4.3,4.3){3} {\line(1,0){4.3}}

\put(137.6,68.8){\line(0,1){4.3}}
\put(137.6,73.1){\line(1,0){4.3}}
\put(141.9,73.1){\line(0,1){8.6}}
\multiput(141.9,81.7)(4.3,4.3){2} {\line(1,0){4.3}}
\multiput(146.2,81.7)(4.3,4.3){2} {\line(0,1){4.3}}
\put(150.5,90.3){\line(1,0){8.6}}
\put(159.1,90.3){\line(0,1){4.3}}
\put(159.1,94.6){\line(1,0){4.3}}

\put(206.4,68.8){\line(0,1){8.6}}
\multiput(206.4,77.4)(4.3,4.3){2} {\line(1,0){4.3}}
\multiput(210.7,77.4)(4.3,4.3){2} {\line(0,1){4.3}}
\put(215,86){\line(1,0){8.6}}
\multiput(223.6,86)(4.3,4.3){2} {\line(0,1){4.3}}
\multiput(223.6,90.3)(4.3,4.3){2} {\line(1,0){4.3}}

\multiput(51.6,43)(68.8,0){4} {\vector(1,-1){20}}
\multiput(154.8,77.4)(68.8,0){2} {\vector(1,-1){20}}
\multiput(26.5,29.73)(68.8,0){4} {$_{\nearrow}$}
\multiput(129.7,64.13)(68.8,0){2} {$_{\nearrow}$}

\end{picture}
    \end{center}

\caption{ Auslander-Reiten quiver of $\mathrm{rep}\hspace{0.1cm}kQ/I$.}\label{figure3.8}

\end{figure}

As an example, Figure \ref{figure3.8} shows the Auslander-Reiten quiver of the Nakayama algebra $\mathscr{A}=kQ/I$ associated to the quiver $Q$ shown in Figure \ref{figure3.4} with $I= \langle x_{3}x_{4},x_{1}x_{2}x_{3} \rangle$.

\vspace{1.5cm}

\section{Cluster Variables Associated to Dyck Paths} \label{section3.3}

In this section, we construct an alphabet associated to Dyck paths. And it is given a formula for cluster variables of cluster algebras associated to Dynkin diagrams of type $\mathbb{A}_{n}$.

\subsection{An Alphabet for Dyck Paths}

\begin{figure}[h!]

\begin{center}
\begin{picture}(169,63)

\multiput(0,19.5)(65,0){3} {\line(1,1){39}}

\put(0,19.5){\line(0,1){13}}
\put(0,32.5){\line(1,0){13}}
\put(13,32.5){\line(0,1){26}}
\put(13,58.5){\line(1,0){26}}
\multiput(13,45.5)(4.3,0){3} {\line(1,0){3.3}} 
\multiput(26,45.5)(0,4.3){3} {\line(0,1){3.3}} 

\put(65,19.5){\line(0,1){26}}
\put(65,45.5){\line(1,0){26}}
\put(91,45.5){\line(0,1){13}}
\put(91,58.5){\line(1,0){13}}
\multiput(65,32.5)(4.3,0){3} {\line(1,0){3.3}} 
\multiput(78,32.5)(0,4.3){3} {\line(0,1){3.3}}

\put(130,19.5){\line(0,1){39}}
\put(130,58.5){\line(1,0){39}}
\multiput(130,32.5)(4.3,0){3} {\line(1,0){3.3}} 
\multiput(130,45.5)(4.3,0){6} {\line(1,0){3.3}} 
\multiput(143,32.5)(0,4.3){6} {\line(0,1){3.3}}
\multiput(156,45.5)(0,4.3){3} {\line(0,1){3.3}}

\put(10,0){(a) $U_{1}^{1}$}
\put(75,0){(b) $U_{2}^{1}$}
\put(140,0){(c) $E_{3}$}

\end{picture}
\end{center}

 \caption{Alphabet $H_{3}$.}\label{figure3.9}

\end{figure}
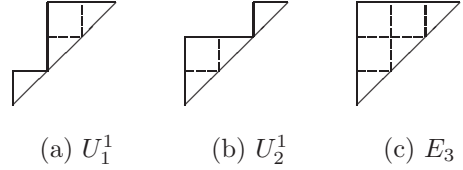

For $n>2$, let $U_{1}^{i}=u_{1} \dots u_{2n}$ and $U_{2}^{i}=u'_{1}\dots u'_{2n}$ be  Dyck paths in $\mathfrak{D}_{2n}$ with the following form:

\begin{equation}
u_{j}=
\begin{cases}
U, & \mbox{ if } 1 \leq j \leq i+1  \text{ or } j=2(i+1)+k \leq 2n,\\
D, & \mbox{ if } i+2 \leq j \leq 2(i+1) \text{ or } j=2(i+1+k)\leq 2n,
\end{cases}
\end{equation}

and 

\begin{equation}
u'_{j}=
\begin{cases}
U, & \mbox{ if } 2i < j \leq i+n  \text{ or }  j=1+2k \leq 2i,\\
D, & \mbox{ if } i+n < j \leq 2n \text{ or } j=2k\leq 2n,
\end{cases}
\end{equation}

for $k >0$ and $ i \leq n-2$. The \textit{alphabet} $H_{n}$ is the union of the set $\lbrace U_{r}^{j} \text{ } |\text{ } r=1,2 \text{ and } 1\leq i \leq n-2 \rbrace$ and the Dyck path with exactly one peak in $\mathfrak{D}_{2n}$ (denoted by $E_{n}$). Figure \ref{figure3.9} shows the alphabet $H_{3}$.\\

Let $\mathscr{C}=\lbrace i_{1}, \dots , i_{k},j_{1}, \dots, j_{m} \rbrace$ be an admissible subchain of \textbf{n-1}. We fix two different relations of concatenation $f_{1}$ and $f_{2}$ over $H_{n}$ such that 

\begin{equation}
f_{1}(V_{i})=
\begin{cases}
E_{n}, & \mbox{ if } V_{i}= E_{n} \text { or } V_{i}=U_{1}^{i},\\
U_{2}^{i+1}, & \mbox{ if }  V_{i}= E_{n} \text { or } V_{i}=U_{1}^{i},\\
U_{1}^{i+1}, & \mbox{ if }  V_{i}=U_{2}^{i},
\end{cases}
\end{equation}

and 

\begin{equation}
f_{2}(V_{i})=
\begin{cases}
E_{n}, & \mbox{ if }  V_{i}=U_{2}^{i},\\
U_{1}^{i+1}, & \mbox{ if }  V_{i}= E_{n} \text { or } V_{i}=U_{1}^{i},\\
U_{2}^{i+1}, & \mbox{ if }  V_{i}=U_{2}^{i}.
\end{cases}
\end{equation}

Then, we take the set of  words $V=V_{1}\dots V_{n-2}$ in $H_{n}^{\ast}$ such that 

\begin{equation}
V_{i}=
\begin{cases}
f_{1}(V_{i-1}) ,& \mbox{ if } i \notin \mathscr{C},  \\
f_{2}(V_{i-1}), & \mbox{ if }  i \in \mathscr{C}-\lbrace1, n-1\rbrace,
\end{cases}
\end{equation}

for $1<i\leq n-2$, $n\geq4$. This set is denoted by $\mathbb{X}_{\mathscr{C}}$, in particular case $\mathbb{X}_{\lbrace 1, 2 \rbrace}=H_{3}$.

\subsection{Dyck Words and Perfect Matchings}

Let $\mathcal{G}=(G_{1}, \dots ,G_{n-1})$ be a snake graph, then we can associate to $\mathcal{G}$ an admissible subchain $\mathscr{C}$ of \textbf{n-1} in the following way:\par\bigskip If $G_{i-1}$, $G_{i}$ and $G_{i+1}$ denote tiles of the following snake graph 

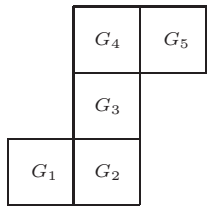
\begin{figure}[h!]

\begin{center}
\begin{picture}(75,75)
\multiput(0,0)(0,25){2} {\line(1,0){50}}
\multiput(25,50)(0,25){2} {\line(1,0){50}}
\put(0,0){\line(0,1){25}}
\multiput(25,0)(25,0){2} {\line(0,1){75}} 
\put(75,50){\line(0,1){25}}

\put(9,12){$_{G_{1}}$}
\put(33,12){$_{G_{2}}$}
\put(33,37){$_{G_{3}}$}
\put(33,62){$_{G_{4}}$}
\put(59,62){$_{G_{5}}$}
 
\end{picture}
\end{center}

 \caption{Snake graph $\mathcal{G}$.}\label{figure3.10}

\end{figure}

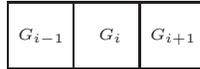
\begin{center}
\begin{picture}(75,25)
\multiput(0,0)(0,25){2} {\line(1,0){75}} 
\multiput(0,0)(25,0){4} {\line(0,1){25}}
\put(4,12){$_{G_{i-1}}$}
\put(34.5,12){$_{G_{i}}$}
\put(54,12){$_{G_{i+1}}$}
\end{picture}
\end{center}

then, $i \in \mathscr{C} $ for $1< i< n-1$. For example, for the snake graph $\mathcal{G}$ shown in Figure \ref{figure3.10}

it holds that the corresponding admissible subchain is given by the identity $\lbrace1,3,5\rbrace= \lbrace i_{1},j_{1},i_{2}\rbrace =\lbrace j_{1},i_{1},j_{2}\rbrace$. By  notation, $\mathcal{G}$ can be written as $\mathcal{G}_{\mathscr{C}}$.\par\bigskip

The following result establishes a relationship between the alphabet $\mathbb{X}_{\mathscr{C}}$ and perfect matchings of snake graphs.

\addtocounter{lema}{8}

\begin{lema}\label{lema3.6}

\textit{Let $\mathscr{C}= \lbrace i_{1}, \dots , i_{k},j_{1}, \dots, j_{m} \rbrace$ be an admissible subchain of \textbf{n-1}. Then, there is a bijective correspondence  between the set $\mathbb{X}_{\mathscr{C}}$ and the perfect matchings of $\mathcal{G}_{\mathscr{C}}$.}

\end{lema} 

\textbf{Proof.} Let $\mathscr{C}$ be an admissible subchain of \textbf{n-1},  $\mathbb{X}_{\mathscr{C}}$ be a set of words, and $\mathcal{G}_{\mathscr{C}}$ be a snake graph associated to $\mathscr{C}$. Assume a numbering over the edges of $\mathcal{G}_{\mathscr{C}}$ in the following way:\par\bigskip 

For boundary edges of  $G_{i}$, we have the following four possibilities  

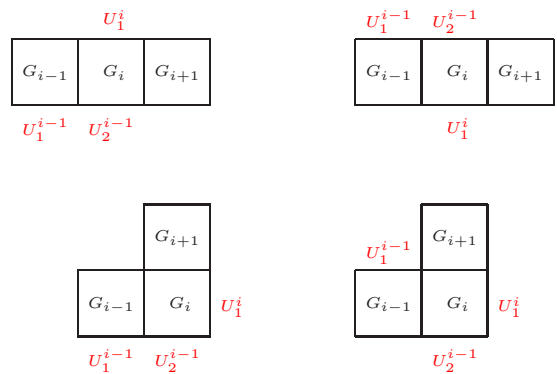
\begin{center} 

\begin{picture}(205,140)

\multiput(25,10)(0,25){2} {\line(1,0){50}}
\multiput(50,10)(25,0){2} {\line(0,1){50}}
\put(25,10){\line(0,1){25}}
\put(50,60){\line(1,0){25}}
\multiput(29,22) (105,0){2}{$_{G_{i-1}}$}
\multiput(59.5,22)(105,0){2}{$_{G_{i}}$}
\multiput(54,47)(105,0){2}{$_{G_{i+1}}$}
\multiput(29,2)(105,41){2}{\textcolor{red}{$_{U_{1}^{i-1}}$}}
\multiput(54,2)(105,0){2}{\textcolor{red}{$_{U_{2}^{i-1}}$}}
\multiput(79,22)(105,0){2}{\textcolor{red}{$_{U_{1}^{i}}$}}

\multiput(130,10)(0,25){2} {\line(1,0){50}}
\multiput(155,10)(25,0){2} {\line(0,1){50}}
\put(130,10){\line(0,1){25}}
\put(155,60){\line(1,0){25}}

\multiput(0,97.5)(0,25){2} {\line(1,0){75}}
\multiput(0,97.5)(25,0){4} {\line(0,1){25}}
\multiput(4,109.5)(130,0){2}{$_{G_{i-1}}$}
\multiput(34.5,109.5)(130,0){2}{$_{G_{i}}$}
\multiput(54,109.5)(130,0){2}{$_{G_{i+1}}$}
\multiput(4,89.5)(130,41){2}{\textcolor{red}{$_{U_{1}^{i-1}}$}}
\multiput(29,89.5)(130,41){2}{\textcolor{red}{$_{U_{2}^{i-1}}$}}
\multiput(34.5,130.5)(130,-41){2}{\textcolor{red}{$_{U_{1}^{i}}$}}

\multiput(130,97.5)(0,25){2} {\line(1,0){75}}
\multiput(130,97.5)(25,0){4} {\line(0,1){25}}

\end{picture}

\end{center}

with $1<i<n-1$ (labeling is given by recurrence). The  other edges are labeled with  the letter $E_{n}$. Now, a perfect matching $P$ of $\mathcal{G}_{\mathscr{C}}$  can be written as a vector $v=(v_{1},\dots, v_{n})$, where each $v_{i}$ corresponds to an edge of $\mathcal{G}_{\mathscr{C}}$ (this vector is unique up to permutation). Define a map $f:\mathbb{X}_{\mathscr{C}}\rightarrow \text{Match}(\mathcal{G}_{\mathscr{C}})$ such that $f(V_{1}\dots V_{n-2})=(E_{n},V_{1},\dots,V_{n-2},E_{n})$. Firstly, we will prove that $f$ is well defined by induction over $n$. To start note that for $n=3$, we have the following three cases: 

\begin{enumerate}

\item [(I)] If $V_{1}=E_{3}$,  it turns out that $F(V_{1})=(E_{3}, E_{3},E_{3})$,  which is given by

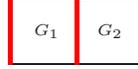
\begin{center}
\begin{picture}(50,30)
\multiput(0,0)(0,25){2} {\line(1,0){50}}
\put(9,12){$_{G_{1}}$}
\put(33,12){$_{G_{2}}$}
\linethickness{1.5pt}
\multiput(0,0)(25,0){3} {\textcolor{red}{\line(0,1){25}}}
\end{picture}

\end{center}

\item [(II)] If $V_{1}=U_{1}^{1}$, it holds that $f(U_{1}^{1})=(E_{3},U_{1}^{1},E_{3})$, which is equal to

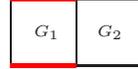
\begin{center}

\begin{picture}(50,30)
\multiput(0,0)(0,25){2} {\line(1,0){50}}
\multiput(0,0)(25,0){3} {\line(0,1){25}}
\put(9,12){$_{G_{1}}$}
\put(33,12){$_{G_{2}}$}
\linethickness{1.5pt}
\multiput(0,0)(0,25){2} {\textcolor{red}{\line(1,0){25}}}
\put(50,0){\textcolor{red}{\line(0,1){25}}}
\end{picture}

\end{center}

\item [(III)] If $V_{1}=U_{2}^{1}$,  then $f(U_{2}^{1})=(E_{3},U_{2}^{1},E_{3})$, which is of the form

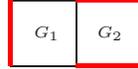
\begin{center}
\begin{picture}(50,30)
\multiput(0,0)(0,25){2} {\line(1,0){50}}
\multiput(0,0)(25,0){3} {\line(0,1){25}}
\put(9,12){$_{G_{1}}$}
\put(33,12){$_{G_{2}}$}
\linethickness{1.5pt}
\multiput(25,0)(0,25){2} {\textcolor{red}{\line(1,0){25}}}
\put(0,0){\textcolor{red}{\line(0,1){25}}}
\end{picture}
\end{center}

\end{enumerate}
Suppose that the result holds for $n=k$. Let $n=k+1$, by hypothesis $(E_{k+1},V_{1}, \dots V_{k})$  are disjoint sets containing all the previous tiles in $\mathcal{G}_{\mathscr{C}}$, then  there are two possibilities for $k$.

\begin{enumerate}

\item [(I)] for $k \in \mathscr{C}- \lbrace 1,k+1 \rbrace$, we have the following conditions:

\begin{enumerate}

\item [(1.1)] If $V_{k-1}=E_{k+1}$, then $f(V_{1} \dots E_{k+1}E_{k+1})=(E_{k+1},V_{1}, \dots, E_{k+1},E_{k+1}, E_{k+1})$ and  $f(V_{1} \dots E_{k+1} U_{2}^{k})=(E_{k+1},V_{1}, \dots, E_{k+1}, U_{2}^{k}, E_{k+1})$, which are given by

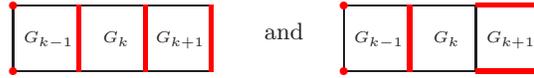
\begin{center}
\begin{picture}(200,30)
\multiput(0,0)(0,25){2} {\line(1,0){75}}
\multiput(0,0)(25,0){4} {\line(0,1){25}}

\multiput(125,0)(0,25){2} {\line(1,0){75}}
\multiput(125,0)(25,0){4} {\line(0,1){25}}

\multiput(4,12)(125,0){2}{$_{G_{k-1}}$}
\multiput(34,12)(125,0){2}{$_{G_{k}}$}
\multiput(54,12)(125,0){2}{$_{G_{k+1}}$}
\put(95,12){and}

\multiput(0,0)(0,25){2}{\textcolor{red}{\circle*{3}}}
\multiput(125,0)(0,25){2}{\textcolor{red}{\circle*{3}}}

\linethickness{1.5pt}
\multiput(25,0)(25,0){3} {\textcolor{red}{\line(0,1){25}}}
\multiput(175,0)(0,25){2} {\textcolor{red}{\line(1,0){25}}}
\put(150,0){\textcolor{red}{\line(0,1){25}}}
\end{picture}

\end{center}

\item [(1.2)] If $V_{k-1}=U_{1}^{k-1}$, then $f(V_{1} \dots U_{1}^{k-1}E_{k+1})=(E_{k+1},V_{1}, \dots, U_{1}^{k-1},E_{k+1}, E_{k+1})$ and $f(V_{1} \dots U_{1}^{k-1}U_{2}^{k})=(E_{k+1},V_{1}, \dots, U_{1}^{k-1}, U_{2}^{k}, E_{k+1})$, which are equal to 

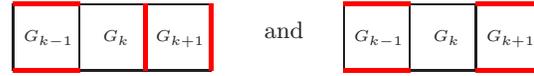
\begin{center}
\begin{picture}(200,30)
\multiput(0,0)(0,25){2} {\line(1,0){75}}
\multiput(0,0)(25,0){4} {\line(0,1){25}}

\multiput(125,0)(0,25){2} {\line(1,0){75}}
\multiput(125,0)(25,0){4} {\line(0,1){25}}

\multiput(4,12)(125,0){2}{$_{G_{k-1}}$}
\multiput(34,12)(125,0){2}{$_{G_{k}}$}
\multiput(54,12)(125,0){2}{$_{G_{k+1}}$}
\put(95,12){and}

\linethickness{1.5pt}
\multiput(0,0)(0,25){2} {\textcolor{red}{\line(1,0){25}}}
\multiput(50,0)(25,0){2} {\textcolor{red}{\line(0,1){25}}}
\multiput(125,0)(0,25){2} {\textcolor{red}{\line(1,0){25}}}
\multiput(175,0)(0,25){2} {\textcolor{red}{\line(1,0){25}}}
\end{picture}
\end{center}

\item [(1.3)] If $V_{k-1}=U_{2}^{k-1}$, then $f(V_{1} \dots U_{2}^{k-1} U_{1}^{k})=(E_{k+1},V_{1}, \dots, U_{2}^{k-1}, U_{1}^{k}, E_{k+1})$   which is of the form

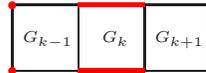
\begin{center}

\begin{picture}(75,30)
\multiput(0,0)(0,25){2} {\line(1,0){75}}
\multiput(0,0)(25,0){4} {\line(0,1){25}}

\put(4,12){$_{G_{k-1}}$}
\put(34,12){$_{G_{k}}$}
\put(54,12){$_{G_{k+1}}$}

\multiput(0,0)(0,25){2}{\textcolor{red}{\circle*{3}}}

\linethickness{1.5pt}
\multiput(25,0)(0,25){2} {\textcolor{red}{\line(1,0){25}}}
\put(75,0){\textcolor{red}{\line(0,1){25}}}
\end{picture}
\end{center}

\end{enumerate}

\item [(II)] for $k  \notin \mathscr{C}$, there are the following cases:

\begin{itemize}

\item [(2.1)] If $V_{k-1}=E_{k+1}$, then $f(V_{1} \dots E_{k+1}U_{1}^{k})=(E_{k+1},V_{1}, \dots, E_{k+1},U_{1}^{k},E_{k+1})$, which is given by

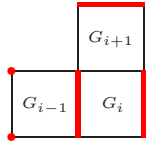
\begin{center}

\begin{picture}(50,55)
\multiput(0,0)(0,25){2} {\line(1,0){50}}
\multiput(25,0)(25,0){2} {\line(0,1){50}}
\put(0,0){\line(0,1){25}}
\put(25,50){\line(1,0){25}}

\put(4,12) {$_{G_{i-1}}$}
\put(34,12){$_{G_{i}}$}
\put(29,37){$_{G_{i+1}}$}

\multiput(0,0)(0,25){2}{\textcolor{red}{\circle*{3}}}

\linethickness{1.5pt}
\multiput(25,0)(25,0){2} {\textcolor{red}{\line(0,1){25}}}
\put(25,50){\textcolor{red}{\line(1,0){25}}}

\end{picture}
\end{center}

\item [(2.2)] If $V_{k-1}=U_{1}^{k-1}$, then $f(V_{1} \dots U_{1}^{k-1}U_{1}^{k})=(E_{k+1},V_{1}, \dots, U_{1}^{k-1},U_{1}^{k},E_{k+1})$, which is equal to

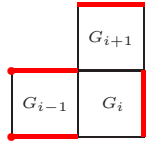
\begin{center}

\begin{picture}(50,55)
\multiput(0,0)(0,25){2} {\line(1,0){50}}
\multiput(25,0)(25,0){2} {\line(0,1){50}}
\put(0,0){\line(0,1){25}}
\put(25,50){\line(1,0){25}}

\put(4,12) {$_{G_{i-1}}$}
\put(34,12){$_{G_{i}}$}
\put(29,37){$_{G_{i+1}}$}

\multiput(0,0)(0,25){2}{\textcolor{red}{\circle*{3}}}

\linethickness{1.5pt}
\multiput(0,0)(0,25){2} {\textcolor{red}{\line(1,0){25}}}
\put(50,0){\textcolor{red}{\line(0,1){25}}}
\put(25,50){\textcolor{red}{\line(1,0){25}}}

\end{picture}

\end{center}

\item [(2.3)] If $V_{k-1}=U_{2}^{k-1}$, then $f(V_{1} \dots U_{2} ^{k-1}E_{k+1})=(E_{k+1},V_{1}, \dots, U_{2}^{k-1}, E_{k+1}, E_{k+1})$  and  $f(V_{1} \dots U_{2}^{k-1}U_{2}^{k})=(E_{k+1},V_{1}, \dots, U_{2}^{k-1}, U_{2}^{k}, E_{k+1})$,  which are of the form

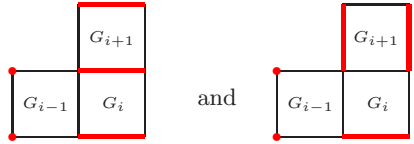
\begin{center}

\begin{picture}(150,55)
\multiput(0,0)(0,25){2} {\line(1,0){50}}
\multiput(25,0)(25,0){2} {\line(0,1){50}}
\put(0,0){\line(0,1){25}}
\put(25,50){\line(1,0){25}}

\multiput(100,0)(0,25){2} {\line(1,0){50}}
\multiput(125,0)(25,0){2} {\line(0,1){50}}
\put(100,0){\line(0,1){25}}
\put(125,50){\line(1,0){25}}

\multiput(4,12)(100,0){2} {$_{G_{i-1}}$}
\multiput(34,12)(100,0){2}{$_{G_{i}}$}
\multiput(29,37)(100,0){2}{$_{G_{i+1}}$}
\put(70,12){and}

\multiput(0,0)(0,25){2}{\textcolor{red}{\circle*{3}}}
\multiput(100,0)(0,25){2}{\textcolor{red}{\circle*{3}}}

\linethickness{1.5pt}
\multiput(25,0)(0,25){3} {\textcolor{red}{\line(1,0){25}}}
\multiput(125,25)(25,0){2} {\textcolor{red}{\line(0,1){25}}}
\put(125,0){\textcolor{red}{\line(1,0){25}}}

\end{picture}

\end{center}

\end{itemize}

\end{enumerate}

Dual arguments prove the result for the other labelings. We also note that by definition map $f$ is injective and  surjective.$\hfill\square$\\

\addtocounter{Nota}{3}

\begin{Nota}

Each perfect matching of $\mathcal{G}_{\mathscr{C}}$ is in correspondence with just only one object of the $\mathbb{A}_{n-1}-$Dyck paths category associated  to the admissible subchain $\mathscr{C}=\lbrace i_{1}, \dots, i_{k}, j_{1}, \dots, j_{m} \rbrace$.

\end{Nota}

For each Dyck path $Y=y_{1} \dots y_{2n}$ with $n-1$ peaks, we construct a family of words $Y\cap\mathbb{X}_{\mathscr{C}}\in H_{n}^{\ast}$ such that;

\begin{equation}
Y \cap \mathbb{X}_{\mathscr{C}}=\lbrace Y \cap V^{z} \text { } | \text { } V^{z} \in  \mathbb{X}_{\mathscr{C}} \rbrace,
\end{equation}

where 

\begin{equation}
Y \cap V^{z}=
\begin{cases}
V^{z} & \mbox{ if there exists } j \text{ such that } y_{j}= v^{z}_{j} \text{ for } 1<j<2n,  \\
E_{n} & \mbox{ otherwise,}
\end{cases}
\end{equation}

with $V^{z}=v^{z}_{1} \dots v^{z}_{2n}$ in $\mathbb{X}_{\mathscr{C}}$. For the set $Y \cap \mathbb{X}_{\mathscr{C}}$, it can be defined a relation $\backsim$ such that

\begin{equation}
Y \cap V^{z_{1}} \backsim Y \cap V^{z_{2}} \text{ if and only if } Y \cap V^{z_{1}}  \text{ and }Y \cap V^{z_{2}} \text{ are the same word. }
\end{equation}

In this case, $\backsim$ is an equivalence relation and $ (Y \cap \mathbb{X}_{\mathscr{C}}) / \backsim $ is denoted by $[Y \cap \mathbb{X}_{\mathscr{C}} ]$. Also, we remind that a Dyck path $Y$ can be written as the word $UWD=Uw_{1}, \dots w_{n-1}D$, where $y_{1}=U$, $y_{2n}=D$ and, $w_{i}=y_{2i}y_{2i+1}$.

\addtocounter{lema}{1}

\begin{lema}\label{lema3.7}

\textit{Let $\mathscr{C}= \lbrace i_{1}, \dots , i_{k},j_{1}, \dots, j_{m} \rbrace$ be an admissible subchain of \textbf{n-1} and let $Y$ a Dyck path of length $2n$ with exactly  $n-1$ peaks. Then, there is a bijective correspondence  between  the set $[Y \cap \mathbb{X}_{\mathscr{C}} ]$ and the set of perfect matchings of the snake graph belonging to $\mathcal{G}_{\mathscr{C}}$ and induced by the words $w_{t}=UD$ in $Y$.}

\end{lema}

\textbf{Proof.} Let $\mathscr{C}$ be an admissible subchain of \textbf{n-1} and $Y=UWD$ be a Dyck path in $S$, then  by Proposition  \ref{prop3.2} there are  $l, r \in \mathbf{Z}_{>0}$ with $1\leq l \leq r \leq n-1$  such that  $w_{t}=UD$ for $l \leq t \leq r$ and $w_{t}=DU$ otherwise. Now, let $\mathcal{G}_{\mathscr{C}^{l, r}}=\mathcal{G}[l, d]$ be a snake graph  belonging  to $\mathcal{G}_{\mathscr{C}}$ induced by $Y$. Define a map $g:[Y \cap \mathbb{X}_{\mathscr{C}} ] \rightarrow \mathrm{Match}(\mathcal{G}_{\mathscr{C}^{l, r}})$ such that:

\begin{enumerate}

\item [(I)] If $1<l\leq r< n-1$, then $g([Y\cap V^{i}])=g(E_{n}\dots E_{n}V^{i}_{l-1}\dots V^{i}_{r}E_{n}\dots E_{n})=(V^{i}_{l-1},\dots ,V^{i}_{r})$.

\item [(II)] If $l=1$ and $1=l \leq r < n-1$, then $g([Y\cap V^{i}])=g(V^{i}_{1}\dots V^{i}_{r}E_{n}\dots E_{n})=(E_{n},V^{i}_{l},\dots ,V^{i}_{r})$.

\item [(III)] If $r=n-1$ and $1< l \leq r=n-1$, then $g([Y\cap V^{i}])=g(E_{n}\dots E_{n}V^{i}_{l-1}\dots V^{i}_{n-2})=(V^{i}_{l-1},\dots ,V^{i}_{n-2},E_{n})$.

\item [(IV)] If $l=1$ and $r=n-1$, then $g=f$.

\end{enumerate}

Since in the four cases $g$ is a restriction of $f$. It follows that $g$ is a bijection as a consequence of Lemma \ref{lema3.6}.$\hfill\square$

\subsection{Cluster Variables Formula Based on Dyck Paths Categories}

In this section, Dyck paths categories are used to give a formula for cluster variables of cluster algebras of Dynkin type $\mathbb{A}_{n}$, to do that, we use the category of Dyck paths associated to an admissible subchain.  We also present a connection between cluster variables of algebras of type $\mathbb{A}_{n-1}$ and Dyck paths with $n-1$ peaks.\par\bigskip

Let $\mathscr{C}=\lbrace i_{1}, \dots i_{k},j_{1}, \dots j_{m}\rbrace$ be an admissible subchain of \textbf{n-1}  and let $Y=UWD$ be a Dyck path in $S$, then we define the monomials

\begin{equation} \label{subsection3.3.3.1}
\eta_{Y}=\prod_{UD=w_{i} \in Y}x_{i},
\end{equation}

and

\begin{equation} \label{subsection3.3.3.2}
X_{V}=\prod_{m \in M_{V}}x_{m},
\end{equation}

with $M_{V}$ being the set of indices $m$ such that

\begin{equation} \label{subsection3.3.3.3}
m=
\begin{cases}
i+1, & \mbox{ if } U_{1}^{i} \in V, \\
i,& \mbox{ if } U_{2}^{i} \in V,\\
0, & \mbox { if } E_{n} \in V,
\end{cases}
\end{equation}

  $V \in [Y \cap \mathbb{X}_{\mathscr{C}} ]$. For this case $x_{0}=1$.\\

The following theorem gives the cluster variable associated to a Dyck path in the set $S$ and its connection with cluster algebras of type $\mathbb{A}_{n-1}$.

 \addtocounter{teor}{4}

\begin{teor}

\textit{Let $\mathscr{C}=\lbrace i_{1}, \dots ,i_{k},j_{1}, \dots ,j_{m}\rbrace$ be an admissible subchain of \textbf{n-1}, $Y=UWD$ a Dyck path with $n-1$ peaks and $M$ the set of  all cluster  variables of a cluster algebra  of type $\mathbb{A}_{n-1}$  with $\lbrace i_{1}, \dots, i_{k} \rbrace$ and $\lbrace  j_{1}, \dots, j_{m} \rbrace$ the sets of sinks and sources, respectively. Then}:

\begin{enumerate}

\item [(i)] \textit{The cluster variable associated  to $Y$ in the category $\mathfrak{C}_{2n}$ is given by} 

\begin{equation}
X_{Y}=(\eta_{Y})^{-1} \Bigg (\sum_{V \in [Y \cap \mathbb{X}_{\mathscr{C}} ]}X_{V} \Bigg ).
\end{equation}

\item  [(ii)]\textit{There exists a bijective correspondence  between Dyck paths with $n-1$ peaks and  the set $M\setminus \text{\textbf{x}}_{0}$ with $\text{\textbf{x}}_{0}$ the initial seed.} 

\end{enumerate}

\end{teor}

\textbf{Proof.} Let $\mathscr{C}=\lbrace i_{1}, \dots ,i_{k},j_{1}, \dots ,j_{m}\rbrace$ be an admissible subchain of \textbf{n-1}, and let $T_{\mathscr{C}}$ be   the triangulation of the polygon with $n+2$ vertices given by $\mathscr{C}$

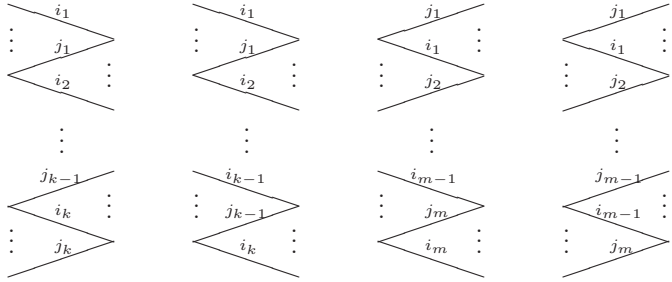
\begin{center}

\begin{picture}(250,110)
\multiput(0,0)(0,26.5){2} {\line(3,1){40}}
\multiput(40,13.3)(0,49.8){2} {\line(-3,1){40}}
\multiput(0,76.3)(70,0){2} {\line(3,1){40}}
\multiput(40,90)(70,0){2} {\line(-3,1){40}}
\multiput(110,26.6)(0,36.5){2} {\line(-3,1){40}}
\multiput(70,13.3)(70,0){2} {\line(3,1){40}}
\multiput(110,0)(70,0){2} {\line(-3,1){40}}
\multiput(180,26.6)(0,49.9){2} {\line(-3,1){40}}
\multiput(140,62.8)(0,27){2} {\line(3,1){40}}
\multiput(210,0)(0,26.5){2} {\line(3,1){40}}
\multiput(250,13.3)(0,63.1){2} {\line(-3,1){40}}
\multiput(210,62.8)(0,27){2} {\line(3,1){40}}

\multiput(0,9)(0,76){2} {$\vdots$}
\multiput(37,22.3)(0,49.5){2} {$\vdots$}
\multiput(70,22.3)(0,62.7){2} {$\vdots$}
\multiput(107,9)(0,62.7){2} {$\vdots$}
\multiput(140,22.3)(0,49.5){2} {$\vdots$}
\multiput(177,9)(0,76){2} {$\vdots$}
\multiput(210,9)(0,62.7){2} {$\vdots$}
\multiput(247,22.3)(0,62.7){2} {$\vdots$}
\multiput(19,47)(70,0){4} {$\vdots$}

\multiput(17.5,100)(70,0){2} {$_{i_{1}}$}
\multiput(17.5,86.7)(70,0){2} {$_{j_{1}}$}
\multiput(17.5,73.4)(70,0){2} {$_{i_{2}}$}
\put(12,38){$_{j_{k-1}}$}
\put(17.5,24.5){$_{i_{k}}$}
\put(17.5,11){$_{j_{k}}$}
\put(82,38){$_{i_{k-1}}$}
\put(82,24.5){$_{j_{k-1}}$}
\put(87.5,11){$_{i_{k}}$}
\multiput(157.5,100)(70,0){2} {$_{j_{1}}$}
\multiput(157.5,86.7)(70,0){2} {$_{i_{1}}$}
\multiput(157.5,73.4)(70,0){2} {$_{j_{2}}$}
\put(152,38){$_{i_{m-1}}$}
\put(157.5,24.5){$_{j_{m}}$}
\put(157.5,11){$_{i_{m}}$}
\put(222,38){$_{j_{m-1}}$}
\put(222,24.5){$_{i_{m-1}}$}
\put(227.5,11){$_{j_{m}}$}


\end{picture}

\end{center}

Let $\alpha_{l, r}$ be a diagonal that  is not in $T_{\mathscr{C}}$  that cuts the diagonals $\alpha_{l},\dots \alpha_{r} \in T_{\mathscr{C}}$. We define a functor $\chi: \mathcal{C}_{T_{\mathscr{C}}}\rightarrow \mathfrak{C}_{2n} $ such that $\chi(\alpha_{l, r})=UW_{l, r}D$, where 

\begin{equation}
w_{j}=
\begin{cases}
UD, & \mbox{ if } l\leq j\leq r, \\
DU, & \mbox{ otherwise, }
\end{cases}
\end{equation} 

and  for  any pivoting elementary move $E:\alpha_{r ,l}\rightarrow \alpha'_{r',l'}$ , $\chi(E)$ is  the elementary shift $F=f_{t_{1}}\circ \dots \circ f_{t_{k}}$ from $UW_{l, r}D$ to $UW_{l',r'}D$.  Theorems \ref{teor1.7.1} and \ref{teor3.1} allow us to establish the following sequence of equivalences: 

\begin{equation} \label{equation3.4}
\mathcal{C}_{T_{\mathscr{C}}}\backsimeq \mathrm{Mod}\hspace{0.1cm}Q_{T_{\mathscr{C}}}\backsimeq \mathfrak{C}_{2n},
\end{equation}

therefore $\chi$ is a categorical equivalence. Thus,

\begin{enumerate}

\item [(i)] Functor $\chi$ and Lemma \ref{lema3.7}, allow to establish that $x_{\gamma}=X_{Y}$.

\item [(ii)] The map $\psi:S \rightarrow M\setminus \text{\textbf{x}}_{0}$ such that $\psi(Y)=X_{Y}$ is a bijection as a consequence of Theorem \ref{teor1.6.1} and the definition of functor $\chi$. We are done. $\hfill\square$\\ 

\end{enumerate}

\begin{figure}[h!]

\begin{center} 

\begin{picture}(276,262.5)

\multiput(77,17.5)(42,0){3} {\line(1,1){35}}

\multiput(105,28.5)(42,0){2} {$\oplus$}

\multiput(0,87.5)(0,70){3}{\multiput(0,0)(154,0){2}{\multiput(0,0)(42,0){3} {\line(1,1){35}}}}

\multiput(28,101.5)(0,70){3}{\multiput(0,0)(154,0){2}{\multiput(0,0)(42,0){2} {$\oplus$}}}

\put(77,17.5){\line(0,1){7}}
\put(77,24.5){\line(1,0){7}}
\put(84,24.5){\line(0,1){28}}
\put(84,52.5){\line(1,0){28}}
\multiput(84,31.5)(2.7,0){3} {\line(1,0){1.7}}
\multiput(84,38.5)(2.9,0){5} {\line(1,0){1.9}}
\multiput(84,45.5)(2.7,0){8} {\line(1,0){1.7}}
\multiput(91,31.5)(0,2.7){8} {\line(0,1){1.7}}
\multiput(98,38.5)(0,2.9){5} {\line(0,1){1.7}}
\multiput(105,45.5)(0,2.7){3} {\line(0,1){1.7}}

\put(119,17.5){\line(0,1){21}}
\put(119,38.5){\line(1,0){21}}
\multiput(140,38.5)(7,7){2} {\line(0,1){7}}
\multiput(140,45.5)(7,7){2} {\line(1,0){7}}
\multiput(119,24.5)(2.7,0){3} {\line(1,0){1.7}}
\multiput(119,31.5)(2.9,0){5} {\line(1,0){1.9}}
\multiput(126,24.5)(0,2.9){5} {\line(0,1){1.9}}
\multiput(133,31.5)(0,2.7){3} {\line(0,1){1.7}}

\put(161,17.5){\line(0,1){28}}
\put(161,45.5){\line(1,0){28}}
\put(189,45.5){\line(0,1){7}}
\put(189,52.5){\line(1,0){7}}
\multiput(161,24.5)(2.7,0){3} {\line(1,0){1.7}}
\multiput(161,31.5)(2.9,0){5} {\line(1,0){1.9}}
\multiput(161,38.5)(2.7,0){8} {\line(1,0){1.7}}
\multiput(168,24.5)(0,2.7){8} {\line(0,1){1.7}}
\multiput(175,31.5)(0,2.9){5} {\line(0,1){1.7}}
\multiput(182,38.5)(0,2.7){3} {\line(0,1){1.7}}


\multiput(0,87.5)(154,0){2}{\put(0,0){\line(0,1){14}}
\put(0,14){\line(1,0){14}}
\multiput(14,14)(7,7){3} {\line(0,1){7}}
\multiput(14,21)(7,7){3} {\line(1,0){7}}
\multiput(0,7)(2.7,0){3} {\line(1,0){1.7}}
\multiput(7,7)(0,2.7){3} {\line(0,1){1.7}}

\multiput(42,0)(7,7){2} {\line(0,1){7}}
\multiput(42,7)(7,7){2} {\line(1,0){7}}
\put(56,14){\line(0,1){21}}
\put(56,35){\line(1,0){21}}
\multiput(56,21)(2.7,0){3} {\line(1,0){1.7}}
\multiput(56,28)(2.9,0){5} {\line(1,0){1.9}}
\multiput(63,21)(0,2.9){5} {\line(0,1){1.7}}
\multiput(70,28)(0,2.7){3} {\line(0,1){1.7
}}}

\put(84,87.5){\line(0,1){35}}
\put(84,122.5){\line(1,0){35}}
\multiput(84,94.5)(2.7,0){3} {\line(1,0){1.7}}
\multiput(84,101.5)(2.9,0){5} {\line(1,0){1.9}}
\multiput(84,108.5)(2.7,0){8} {\line(1,0){1.7}}
\multiput(84,115.5)(2.9,0){10} {\line(1,0){1.9}}
\multiput(91,94.5)(0,2.9){10} {\line(0,1){1.9}}
\multiput(98,101.5)(0,2.7){8} {\line(0,1){1.7}}
\multiput(105,108.5)(0,2.9){5} {\line(0,1){1.7}}
\multiput(112,115.5)(0,2.7){3} {\line(0,1){1.7}}


\multiput(0,157.5)(196,0){2}{\put(0,0){\line(0,1){35}}
\put(0,35){\line(1,0){35}}
\multiput(0,7)(2.7,0){3} {\line(1,0){1.7}}
\multiput(0,14)(2.9,0){5} {\line(1,0){1.9}}
\multiput(0,21)(2.7,0){8} {\line(1,0){1.7}}
\multiput(0,28)(2.9,0){10} {\line(1,0){1.9}}
\multiput(7,7)(0,2.9){10} {\line(0,1){1.9}}
\multiput(14,14)(0,2.7){8} {\line(0,1){1.7}}
\multiput(21,21)(0,2.9){5} {\line(0,1){1.7}}
\multiput(28,28)(0,2.7){3} {\line(0,1){1.7}}
}

\multiput(0,157.5)(154,70){2}{\multiput(42,0)(7,7){2} {\line(0,1){7}}
\multiput(42,7)(7,7){2} {\line(1,0){7}}
\put(56,14){\line(0,1){21}}
\put(56,35){\line(1,0){21}}
\multiput(56,21)(2.7,0){3} {\line(1,0){1.7}}
\multiput(56,28)(2.9,0){5} {\line(1,0){1.9}}
\multiput(63,21)(0,2.9){5} {\line(0,1){1.7}}
\multiput(70,28)(0,2.7){3} {\line(0,1){1.7
}}}

\multiput(84,157.5)(154,-70){2}{\multiput(0,0)(7,7){3} {\line(0,1){7}}
\multiput(0,7)(7,7){3} {\line(1,0){7}}
\put(21,21){\line(0,1){14}}
\put(21,35){\line(1,0){14}}
\multiput(21,28)(2.7,0){3} {\line(1,0){1.7}}
\multiput(28,28)(0,2.7){3} {\line(0,1){1.7}}
}

\put(154,157.5){\line(0,1){14}}
\put(154,171.5){\line(1,0){14}}
\multiput(168,171.5)(7,7){3} {\line(0,1){7}}
\multiput(168,178.5)(7,7){3} {\line(1,0){7}}
\multiput(154,164.5)(2.7,0){3} {\line(1,0){1.7}}
\multiput(161,164.5)(0,2.7){3} {\line(0,1){1.7}}


\multiput(0,227.5)(154,0){2}{\put(0,0){\line(0,1){35}}
\put(0,35){\line(1,0){35}}
\multiput(0,7)(2.7,0){3} {\line(1,0){1.7}}
\multiput(0,14)(2.9,0){5} {\line(1,0){1.9}}
\multiput(0,21)(2.7,0){8} {\line(1,0){1.7}}
\multiput(0,28)(2.9,0){10} {\line(1,0){1.9}}
\multiput(7,7)(0,2.9){10} {\line(0,1){1.9}}
\multiput(14,14)(0,2.7){8} {\line(0,1){1.7}}
\multiput(21,21)(0,2.9){5} {\line(0,1){1.7}}
\multiput(28,28)(0,2.7){3} {\line(0,1){1.7}}}

\multiput(42,227.5)(196,0){2}{\put(0,0){\line(0,1){35}}
\put(0,35){\line(1,0){35}}
\multiput(0,7)(2.7,0){3} {\line(1,0){1.7}}
\multiput(0,14)(2.9,0){5} {\line(1,0){1.9}}
\multiput(0,21)(2.7,0){8} {\line(1,0){1.7}}
\multiput(0,28)(2.9,0){10} {\line(1,0){1.9}}
\multiput(7,7)(0,2.9){10} {\line(0,1){1.9}}
\multiput(14,14)(0,2.7){8} {\line(0,1){1.7}}
\multiput(21,21)(0,2.9){5} {\line(0,1){1.7}}
\multiput(28,28)(0,2.7){3} {\line(0,1){1.7}}}

\multiput(84,227.5)(154,-70){2}{\put(0,0){\line(0,1){28}}
\put(0,28){\line(1,0){28}}
\put(28,28){\line(0,1){7}}
\put(28,35){\line(1,0){7}}
\multiput(0,7)(2.7,0){3} {\line(1,0){1.7}}
\multiput(0,14)(2.9,0){5} {\line(1,0){1.9}}
\multiput(0,21)(2.7,0){8} {\line(1,0){1.7}}
\multiput(7,7)(0,2.7){8} {\line(0,1){1.7}}
\multiput(14,14)(0,2.9){5} {\line(0,1){1.7}}
\multiput(21,21)(0,2.7){3} {\line(0,1){1.7}}
}

\put(25,210){$_{\text{(a) } E_{5}\oplus E_{5}\oplus U_{1}^{3}}$}
\put(179,210){$_{\text{(b) } E_{5}\oplus U_{2}^{2}\oplus E_{5}}$}
\put(25,140){$_{\text{(c) } E_{5}\oplus U_{2}^{2}\oplus U_{2}^{3}}$}
\put(179,140){$_{\text{(d) } U_{1}^{1}\oplus E_{5}\oplus U_{1}^{3}}$}
\put(25,70){$_{\text{(e) } U_{1}^{1}\oplus U_{2}^{2}\oplus E_{5}}$}
\put(179,70){$_{\text{(f) } U_{1}^{1}\oplus U_{2}^{2}\oplus U_{2}^{3}}$}
\put(102,0){$_{\text{(h) } U_{2}^{1}\oplus U_{1}^{2}\oplus U_{1}^{3}}$}

\end{picture}

\end{center}

\caption{Objects in $\mathfrak{C}_{10}$.}\label{figure3.11}

\end{figure}
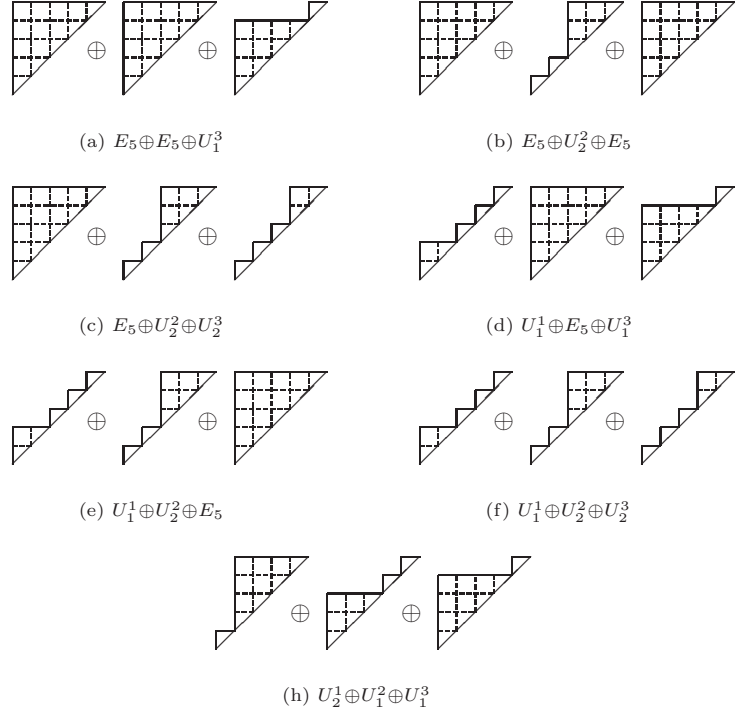

For instance, let $\mathscr{C}=\lbrace j_{1}=1, i_{1}=2, j_{2}=4\rbrace$ be an admissible subchain  of \textbf{4}, the set $\mathbb{X}_{\mathscr{C}}$ is in correspondence with the objects of $\mathfrak{C}_{10}$ shown in Figure \ref{figure3.11}.\\

Then, for $Y=UDUUDUDDUD$, we define the set $Y\cap\mathbb{X}_{\mathscr{C}}$ such that 

\begin{equation}
[Y \cap \mathbb{X}_{\mathscr{C}}]=\lbrace E_{5}E_{5}U_{1}^{3}, E_{5}U_{2}^{2}E_{5},U_{2}^{1}U_{1}^{2}U_{1}^{3} \rbrace.
\end{equation}

Thus, identities (\ref{subsection3.3.3.1}), (\ref{subsection3.3.3.2}) and (\ref{subsection3.3.3.3}) define the polynomials

\begin{equation}
\eta_{Y}=x_{2}x_{3}, \text{ } X_{E_{5}E_{5}U_{1}^{3}}=x_{0}x_{0}x_{4}, \text{ } X_{E_{5}U_{2}^{2}E_{5}}=x_{0}x_{2}x_{0}, \text{ } X_{U_{2}^{1}U_{1}^{2}U_{1}^{3}}=x_{3}x_{1}x_{4},
\end{equation}

 therefore, the cluster variable associated to the Dyck path $Y$ is given by the expression

 \begin{equation}
 X_{Y}=\frac{x_{4}+x_{2}+x_{3}x_{1}x_{4}}{x_{2}x_{3}}.
 \end{equation}


\end{document}